\documentclass[11pt]{amsart}

\usepackage{amscd,amssymb}
\usepackage{amsthm,amsmath,amssymb}
\usepackage[matrix,arrow]{xy}

\makeatletter\@addtoreset{equation}{section} \makeatother

\newcommand{\xref}[1]{{\rm \ref{#1}}}
\newcommand{\Supp}{\operatorname{Supp}}

\newcommand{\Sing}{\operatorname{Sing}}
\newcommand{\Gal}{\operatorname{Gal}}
\newcommand{\Proj}{\operatorname{Proj}}

\newcommand{\Pic}{\operatorname{Pic}}
\newcommand{\mult}{\operatorname{mult}}
\newcommand{\id}{\operatorname{id}}
\newcommand{\Aut}{\operatorname{Aut}}
\newcommand{\Cl}{\operatorname{Cl}}
\newcommand{\cha}{\operatorname{char}}

\newcommand{\Bir}{\operatorname{Bir}}

\newcommand{\codim}{\operatorname{codim}}

\newcommand{\rank}{\operatorname{rank}}
\renewcommand{\emptyset}{\varnothing}

\newcommand{\NN}{{\mathbb N}}
\newcommand{\ZZ}{{\mathbb Z}}
\newcommand{\CC}{{\mathbb C}}

\newcommand{\QQ}{{\mathbb Q}}
\newcommand{\PP}{{\mathbb P}}
\newcommand{\KK}{{\mathbb K}}
\newcommand{\RR}{{\mathbb R}}
\newcommand{\FF}{{\mathbb F}}
\newcommand{\CCSS}{{\mathbb {CS}}}
\newcommand{\LLCCSS}{{\mathbb {LCS}}}

\newcommand{\HHH}{\mathcal H}

\newcommand{\MMM}{\mathcal M}
\newcommand{\OOO}{\mathcal O}
\newcommand{\DDD}{\mathcal D}

\newcommand{\III}{\mathcal I}
\newcommand{\LLL}{\mathcal L}
\newcommand{\NNN}{\mathcal N}
\newcommand{\TTT}{\mathcal T}

\newtheorem{theorem}[equation]{Theorem}

\newtheorem{proposition}[equation]{Proposition}
\newtheorem{lemma}[equation]{Lemma}
\newtheorem{corollary}[equation]{Corollary}

\newtheorem{theoremA}{Theorem A}
\newtheorem{theoremB}{Theorem B}
\newtheorem{theoremC}{Theorem C}
\newtheorem{theoremD}{Theorem D}
\newtheorem{theoremE}{Theorem E}
\newtheorem{theoremF}{Theorem F}
\newtheorem{corollaryA}{Corollary A}
\newtheorem{corollaryB}{Corollary B}
\newtheorem{corollaryC}{Corollary C1}
\newtheorem{corollaryCC}{Corollary C2}

\newtheorem{constructionA}{Construction A}
\newtheorem{constructionB}{Construction B}

\theoremstyle{definition}
\newtheorem{example}[equation]{Example}
\newtheorem{definition}[equation]{Definition}

\theoremstyle{remark}
\newtheorem{remark}[equation]{Remark}

\author{Ivan Cheltsov and Jihun Park}
\title{Sextic double solids}

\address{\emph{Ivan Cheltsov}\newline
\textnormal{School of Mathematics, The University of Edinburgh,
Mayfield Road, Edinburgh EH9 3JZ UK;
 \texttt{cheltsov@yahoo.com}}}
\address{
\emph{Jihun Park}\newline \textnormal{Department of Mathematics,
POSTECH, Pohang, Kyungbuk 790-784, Korea;
\texttt{wlog@postech.ac.kr}}}
\begin{document}

\begin{abstract}
We study properties of double covers of $\PP^3$ ramified along
nodal sextic surfaces such as non-rationality, $\QQ$-factoriality,
potential density, and elliptic fibration structures. We also
consider some relevant problems over fields of positive
characteristic.

 \vspace{5mm}

\end{abstract}

\maketitle 


 All varieties are assumed to be projective,
normal, and defined over the field $\CC$ unless otherwise stated.

\section{Introduction.}
For a given variety, it is one of substantial questions whether it
is rational or not. Global holomorphic differential forms are
natural birational invariants of smooth algebraic varieties which
solve the rationality problem for algebraic curves and surfaces
(see \cite{Za58}). However, these birational invariants are not
sensitive enough to tell whether a given higher dimensional
algebraic variety is non-rational. There are only four known
methods to prove the non-rationality of a higher dimensional
algebraic variety (see \cite{Is97}).

The non-rationality of a smooth quartic $3$-fold was proved in
\cite{IsMa71} using the group of birational automorphisms as a
birational invariant. The non-rationality of a smooth cubic
$3$-fold was proved in \cite{ClGr72} through the study of its
intermediate Jacobian. Birational invariance of the torsion
subgroups of the 3rd integral cohomology  groups  were used in
\cite{ArMu72} to prove the non-rationality of some unirational
varieties. The non-rationality of a wide class of rationally
connected varieties was proved in \cite{Ko95} via reductions into
fields of positive characteristic (see \cite{ChWo04}, \cite{Ko96},
and \cite{Ko00}). Meanwhile, the method of intermediate Jacobians
works only in $3$-folds. In most of interesting cases, the 3rd
integral cohomology groups have no torsion. The method of the
paper \cite{Ko95} works in every dimension, but its direct
application gives the non-rationality just for a very general
element of a given family. Even though the method of \cite{IsMa71}
works in every dimension, the area of its application is not so
broad.

For this paper we mainly use the method that has evolved out of
\cite{IsMa71}. The most significant concept in the method is the
birational super-rigidity that was implicitly introduced in
\cite{IsMa71}.
\begin{definition}
\label{definition-of-super-rigidity} A terminal $\QQ$-factorial
Fano variety $V$ with $\Pic(V)\cong \ZZ$ is bi\-ra\-ti\-o\-nal\-ly
super-rigid if the following three conditions hold:
\begin{enumerate}
\item the variety $V$ cannot be birationally transformed into a
fibration\footnote{For every fibration $\tau:Y\to Z$, we  assume
that $\dim(Y)>\dim(Z)\ne 0$ and $\tau_{*}(\OOO_{Y})=\OOO_{Z}$.}
whose general enough fiber is a
smooth variety of Kodaira dimension $-\infty$; %
\item the variety $V$ cannot be birationally transformed into
another terminal $\QQ$-factorial  Fano variety with Picard group
$\ZZ$ that is not biregular to $V$; %
\item $\Bir(V)=\Aut(V)$.
\end{enumerate}
\end{definition}
Implicitly the paper \cite{IsMa71} proved that all the smooth
quartic $3$-folds in $\PP^4$ are birationally super-rigid.
Moreover, some Fano $3$-folds with non-trivial group of birational
automorphisms  were also handled by the technique of
\cite{IsMa71}, which gave the following weakened version of the
birational super-rigidity:

\begin{definition}
\label{definition-of-rigidity} A terminal $\QQ$-factorial Fano
variety $V$ with $\Pic(V)\cong \ZZ$ is called
bi\-ra\-ti\-o\-nal\-ly rigid if the first two conditions of
Definition~\xref{definition-of-super-rigidity} are satisfied.
\end{definition}
It is clear that the birational rigidity implies the
non-rationality. Initially the technique of \cite{IsMa71} was
applied only to smooth varieties such as quartic $3$-folds,
quintic $4$-folds, certain complete intersections, double spaces
and so on, but later, to singular varieties in   \cite{CoMe02},
\cite{Gr88a}, \cite{Gr88b}, \cite{Gr00}, \cite{Me03},
\cite{Pu88b}, \cite{Pu97}, and \cite{Pu03}. Moreover, similar
results were proved for many higher-dimensional conic bundles (see
\cite{Sar80} and \cite{Sar82}) and del Pezzo fibrations (see
\cite{Pu98b}). Recently, Shokurov's connectedness principle in
\cite{Sh92} shed a new light on the birational rigidity, which
simplified the proofs of old results and helped to obtain new
results (see
 \cite{Ch00b}, \cite{Ch03b}, \cite{Co00}, \cite{dFEM03},
 \cite{Pu98a}, and \cite{Pu02}).

A quartic $3$-fold with a single simple double point is not
birationally super-rigid, but it is proved in \cite{Pu88b} to be
birationally rigid (for a simple proof, see \cite{Co00}). However,
a quartic $3$-fold with one non-simple double point may not
necessarily be birationally rigid as shown in \cite{CoMe02}. On
the other hand, $\QQ$-factorial quartic $3$-folds with only simple
double points are birationally rigid (see \cite{Me03}).

Double covers of $\PP^3$ with at most simple double points,
so-called, double solids, were studied in \cite{Cl83} with a
special regard to quartic double solids, \emph{i.e.}, double
covers of $\PP^3$ ramified along quartic nodal surfaces. It is
natural to ask whether a double solid is rational or not. We can
immediately see that all double solids are non-rational when their
ramification surfaces are of degree greater than six. However, if
the ramification surfaces have lower degree, then the problem is
not simple.

Smooth quartic double solids are known to be non-rational (see
\cite{Cl91}, \cite{Le84}, \cite{Ti80a}, \cite{Ti80b}, \cite{Ti82},
\cite{Ti86}, and \cite{We81}), but singular ones can be
birationally transformed into conic bundles. Quartic double solids
cannot have more than $16$ simple double points (see \cite{Be79a},
\cite{En99}, \cite{Ku}, \cite{Ni75}, and \cite{Ro86}) and in the
case of one simple double point they are non-rational as well (see
\cite{Be79b} and \cite{Tu78}). There are non-$\QQ$-factorial
quartic double solids with six simple double points that  can be
birationally transformed into smooth cubic $3$-folds (see
\cite{Kr00}) and therefore are not rational due to \cite{ClGr72}.
On the other hand, some quartic double solids with seven simple
double points are rational (see \cite{Kr00}). In general, the
rationality question of singular quartic double solids can be very
subtle and must be handled through the technique of intermediate
Jacobians (see \cite{Be79b}, \cite{Sh82}, and \cite{Sh83}).

In the present paper we will consider the remaining case -- the
non-rationality question of sextic double solids, \emph{i.e.},
double covers of $\PP^3$ ramified along sextic nodal surfaces. To
generate various examples of sextic double solids, we note that a
double cover $\pi:X\to \PP^3$ ramified along a sextic surface
$S\subset \PP^3$ can be considered as a hypersurface
$$
u^{2}=f_{6}(x,y,z,w)
$$
of degree $6$ in the weighted projective space $\PP(1,1,1,1,3)$,
where $x$, $y$, $z$, and $w$ are homogeneous coordinates of weight
$1$, $u$ is a homogeneous coordinate of weight $3$, and $f_{6}$ is
a homogeneous polynomial of degree $6$.

A smooth sextic double solid is proved to be birationally
super-rigid in \cite{Is80}. Moreover, a smooth double space of
dimension $n\geq 3$ was considered in \cite{Pu88a}. The birational
super-rigidity of a double cover of $\PP^{3}$ ramified along a
sextic with one simple double point was proved in \cite{Pu97}.  To
complete the study in this direction, one needs to prove the
following:
\begin{theoremA}
\label{theorem-main} Let $\pi:X\to \PP^3$ be a $\QQ$-factorial
double cover ramified along a  sextic nodal  surface $S\subset
\PP^3$. Then $X$ is birationally super-rigid.
\end{theoremA}
As an immediate consequence,  we obtain:
\begin{corollaryA}
\label{corollary-main} Every $\QQ$-factorial double cover of
$\PP^3$ ramified along a sextic nodal surface  is non-rational and
not birationally isomorphic to a conic bundle.
\end{corollaryA}

\begin{remark}
\label{remark-non-closed-fields} Our proof of Theorem~A does not
require the base field to be algebraically closed. Therefore, the
statement of Theorems~A is valid over an arbitrary field of
characteristic zero.
\end{remark}
One can try to prove the non-rationality of a sextic double solid
using the technique of intermediate Jacobians (see \cite{Be79b},
\cite{Sh82}, and \cite{Sh83}), but it seems to be very hard and
still undone even in the smooth case (see \cite{CeVe86}) except
for the non-rationality of a sufficiently general smooth sextic
double solid  via a degeneration technique (see \cite{Be79b},
\cite{Cl82}, and  \cite{Tu78}).

 It is worth while to put emphasis on the $\QQ$-factoriality condition of Theorem~A. Indeed,
  rational sextic double solids do exist if we drop the
$\QQ$-factoriality condition.

\begin{example}
\label{example-of-Barth} Let $X$ be the double cover of $\PP^3$
ramified in the Barth sextic  (see \cite{Ba96}) given by the
equation
$$
4(\tau^2x^2-y^2)(\tau^2y^2-z^2)(\tau^2z^2-x^2)-w^2(1+2\tau)(x^2+y^2+z^2-w^{2})^2=0
$$
in $\Proj(\CC[x,y,z,w])$, where $\tau={\frac {1+\sqrt{5}} {2}}$.
Then $X$ has only simple double points and the number of singular
points is $65$. Moreover, there is a determinantal quartic
$3$-fold $V\subset\PP^4$ with $42$ simple double points such that
the diagram
\[\xymatrix{
&V\ar@{-->}[d]_{\rho}\ar@{^{(}->}[rr]^{f}&&\PP^4\ar@{-->}[d]^{\gamma}&\\%
&X\ar@{->}[rr]_{\pi}&&\PP^{3}&}\] %
commutes (see  \cite{En99} and \cite{Pet98}), where $\rho$ is a
birational map and $\gamma$ is the projection from one simple
double point of the quartic $V$. Therefore, the double cover $X$
is rational because determinantal quartics are rational (see
\cite{Me03} and \cite{Pet98}). In particular, $X$ is not
$\QQ$-factorial by Theorem~A. Indeed, one can show that
$\Pic(X)\cong\ZZ$ and $\Cl(X)\cong \ZZ^{14}$ (see Example 3.7 in
\cite{En99}).
\end{example}

 A point $p$ on a double
cover $\pi:X\to \PP^{3}$ ramified along a sextic surface $S$ is a
simple double point on $X$ if and only if the point $\pi(p)$ is a
simple double point on $S$. Sextic surfaces cannot have more than
$65$ simple double points (see \cite{Bas06}, \cite{JaRu97}, and
\cite{Wa98}). Furthermore, for each positive integer $m$ not
exceeding $65$ there is a sextic surface with $m$ simple double
points (see \cite{Ba96}, \cite{CaCe82}, and \cite{St78}), but in
many cases it is not clear whether the corresponding double cover
is $\QQ$-factorial or not (see \cite{Cl83}, \cite{Cy99}, and
\cite{En99}).

Example~\ref{example-of-Barth} shows that the $\QQ$-factoriality
condition is crucial for Theorem~A. Accordingly, it is worth our
while to study the $\QQ$-factoriality of sextic double solids.

A variety $X$ is called $\QQ$-factorial if a multiple of each Weil
divisor on the variety $X$ is a Cartier divisor. The
$\QQ$-factoriality depends on both local types of singularities
and their global position (see \cite{CiGe03}, \cite{Cl83},  and
\cite{Me03}). Moreover, the $\QQ$-factoriality of the variety $X$
depends on the field of definition of the variety $X$ as well.
When $X$ is a Fano $3$-fold with mild singularities and defined
over $\CC$, the global topological condition
$$
\operatorname{rank}(H^{2}(X,\ZZ))=\operatorname{rank}(H_{4}(X,\ZZ))
$$
is equivalent to  the $\QQ$-factoriality. The following three
examples are inspired by \cite{Bar84}, \cite{Kr00},and
\cite{Me03}.

\begin{example}
\label{example-on-Q-factoriality} Let $\pi:X\to \PP^3$ be the
double cover ramified along a sextic $S$ and  given by
$$
u^{2}+g^{2}_{3}(x,y,z,w)=h_{1}(x,y,z,w)f_{5}(x,y,z,w)\subset
\PP(1,1,1,1,3),
$$
where $g_{3}$, $h_{1}$, and $f_{5}$ are sufficiently general
polynomials defined over $\RR$ of degree $3$, $1$, and $5$,
respectively; $x$, $y$, $z$, $w$ are homogeneous coordinates of
weight $1$; $u$ is a homogeneous coordinate of weight $3$. Then
the double cover $X$ is not $\QQ$-factorial over $\CC$ because the
divisor $h_{1}=0$ splits into two non-$\QQ$-Cartier divisors
conjugated by $\Gal(\CC\slash \RR)$ and given by the equation
$$
(u+\sqrt{-1}g_{3}(x,y,z,w))(u-\sqrt{-1}g_{3}(x,y,z,w))=0.
$$

The sextic surface $S\subset \Proj(\CC[x,y,z,w])$ has $15$ simple
double points at the intersection points of the three surfaces
$$
\{h_{1}(x,y,z,w)=0\}\cap\{g_{3}(x,y,z,w)=0\}\cap\{f_{5}(x,y,z,w)=0\},
$$
which gives $15$ simple double points of $X$. Introducing two new
variables  $s$ and $t$ of weight $2$ defined by
$$
\left\{\aligned
&s={\frac{u+\sqrt{-1}g_{3}(x,y,z,w)} {h_{1}(x,y,z,w)}}={\frac{f_{5}(x,y,z,w)} {u-\sqrt{-1}g_{3}(x,y,z,w)}}\\
&t={\frac{u-\sqrt{-1}g_{3}(x,y,z,w)} {h_{1}(x,y,z,w)}}={\frac{f_{5}(x,y,z,w)} {u+\sqrt{-1}g_{3}(x,y,z,w)}}\\
\endaligned
\right.
$$
we can unproject $X\subset \PP(1,1,1,1,3)$ in the sense of
\cite{Re00} into two complete intersections
$$
\left\{\aligned &V_{s}=\left\{\aligned
&s h_{1}(x,y,z,w)=u+\sqrt{-1}g_{3}(x,y,z,w)\\
&s (u-\sqrt{-1}g_{3}(x,y,z,w))=f_{5}(x,y,z,w)\\
\endaligned
\right \}\subset \PP(1,1,1,1,3,2)\\
&V_{t}=\left\{\aligned
&t h_{1}(x,y,z,w)=u-\sqrt{-1}g_{3}(x,y,z,w)\\
&t (u+\sqrt{-1}g_{3}(x,y,z,w))=f_{5}(x,y,z,w)\\
\endaligned
\right\}\subset \PP(1,1,1,1,3,2),\\
\endaligned
\right.
$$
respectively, which are not defined over $\RR$. Eliminating
variable $u$, we get
$$
\left\{\aligned &V_{s}=\{s^2h_{1}-2\sqrt{-1}s g_{3}-f_{5}=0\}\subset \PP(1,1,1,1,2)\\
&V_{t}=\{t^2h_{1}+2\sqrt{-1}t g_{3}-f_{5}=0\}\subset \PP(1,1,1,1,2)\\
\endaligned
\right.
$$
and for the unprojections $\rho_{s}:X\dasharrow V_{s}$ and
$\rho_{t}:X\dasharrow V_{t}$ we obtain a commutative diagram
\[ \xymatrix{
&Y_s\ar[dl]_{\psi_s}  \ar[dr]^{\phi_s} && Y_t \ar[dl]_{\phi_t}\ar[dr]^{\psi_t}&\\
V_s && X\ar@{-->}[rr]^{\rho_t}\ar@{-->}[ll]_{\rho_s}&&V_t}\] with
birational morphisms $\phi_{s}$, $\psi_{s}$, $\phi_{t}$, and
$\psi_{t}$ such that $\psi_{s}$ and $\psi_{t}$ are extremal
contractions in the sense of \cite{Co95}, while $\phi_{s}$ and
$\phi_{t}$ are flopping contractions. Both the weighted
hypersurfaces $V_{s}$ and $V_{t}$ are quasi-smooth (see
\cite{IF00}) and $\QQ$-factorial with Picard groups $\ZZ$
(Lemma~3.5 in \cite{CPR}, Lemma~3.2.2 in \cite{Do82},
Th\'eor\'eme~3.13 of Exp. XI in \cite{Gro65},  see also
\cite{CalLy94}). Moreover, $V_{s}$ and $V_{t}$ are projectively
isomorphic in $\PP(1,1,1,1,2)$ by the action of
$\Gal(\CC\slash\RR)\cong\ZZ_{2}$. In particular,
$$
\Pic(Y_{s})\cong\Pic(Y_{t})\cong \ZZ\oplus\ZZ;
$$
$Y_{s}$ and $Y_{t}$ are $\QQ$-factorial; $\Cl(X)=\ZZ\oplus\ZZ$.
However, the $\Gal(\CC\slash\RR)$-invariant part of the group
$\Cl(X)$ is $\ZZ$. Thus the $3$-fold $X$ is $\QQ$-factorial over
$\RR$. It is therefore birationally super-rigid and non-rational
over $\RR$ by Theorem~A. It  is also not rational over $\CC$
because $V_{s}\cong V_{t}$ is birationally rigid (see \cite{CPR}).
Moreover, the involution of $X$ interchanging fibers of $\pi$
induces a non-biregular involution $\tau\in \Bir(V_{s})$ which is
regularized by $\rho_{s}$, \emph{i.e.}, the self-map
$\rho_s^{-1}\circ\tau\circ\rho_s:X\to X$ is biregular (see
\cite{Ch04c}).
\end{example}

\begin{example}
\label{example-on-quartic-with-double-point} Let $V\subset \PP^4$
be a quartic $3$-fold with one simple double point $o$. Then the
quartic $V$ is $\QQ$-factorial and $\Pic(V)\cong \ZZ$. In fact,
$V$ can be given by the equation
$$
t^{2}f_{2}(x,y,z,w)+tf_{3}(x,y,z,w)+f_{4}(x,y,z,w)=0\subset
\PP^4=\Proj(\CC[x,y,z,w,t]).
$$
Here, the point $o$ is located at $[0:0:0:0:1]$. It is well known
that the quartic $3$-fold $V$ is birationally rigid and hence
non-rational (see \cite{Co00}, \cite{Me03}, and \cite{Pu88b}).
However, the quartic $V$ is not birationally super-rigid because
$\Bir(V)\ne\Aut(V)$. Indeed, the projection $\phi:V\dasharrow
\PP^3$ from the point $o$ has degree $2$ at a generic point of $V$
and induces a non-biregular involution $\tau\in \Bir(V)$.

Let $f:Y\to V$ be the blow up at the point $o$. Then the linear
system $|-nK_{Y}|$ is free for some natural number  $n\gg 0$ and
gives a birational morphism $g=\phi_{|-nK_{Y}|}:Y\to X$
contracting every curve $C_i\subset Y$ such that $f(C_i)$ is a
line on $V$ passing through the point $o$.  We then obtain the
double cover $\pi:X\to\PP^3$ ramified along the sextic surface
$S\subset\PP^3$ given by the equation
$$
f^2_{3}(x,y,z,w)-4f_{2}(x,y,z,w)f_{4}(x,y,z,w)=0.
$$
The variety $X$, \emph{a priori}, has canonical Gorenstein
singularities.

 We suppose that $V$ is general enough. Each line
$f(C_{i})$ then corresponds to an intersection point  of three
surfaces
$$
\{f_{2}(x,y,z,w)=0\}\cap\{f_{3}(x,y,z,w)=0\}\cap\{f_{4}(x,y,z,w)=0\}$$
in $\PP^3=\Proj(\CC[x,y,z,w])$ which gives $24$ different smooth
rational curves $C_{1}, C_2, \cdots, C_{24}$ on $Y$. For each
curve $C_i$ we have
$$
\NNN_{Y\slash C_{i}}\cong \OOO_{C_{i}}(-1)\oplus \OOO_{C_{i}}(-1)
$$
and hence the morphism $g$ is a  standard flopping contraction
which maps every curve $C_{i}$ to a simple double point of the
$3$-fold $X$. In particular, the sextic $S\subset\PP^3$ has
exactly $24$ simple double points. However, the $3$-fold $X$ is
not $\QQ$-factorial and $\Cl(X)=\ZZ\oplus\ZZ$.

Put $\rho:=g\circ f^{-1}$. Then the involution
$\gamma=\rho\circ\tau\circ\rho^{-1}$ is biregular on $X$ and
interchanges the fibers of the double cover $\pi$. Thus the map
$\rho$ is a regularization of the birational non-biregular
involution $\tau$ in the sense of \cite{Ch04c}, while the
commutative diagram
\[ \xymatrix{
&Y\ar[dl]_{f}  \ar[dr]^{g} &&&& Y\ar[dl]_{g}\ar[dr]^{f}&\\
V\ar@{-->}[rr]^{\rho}&&X\ar@{->}[rr]_{\gamma}&&X&&\ar@{-->}[ll]_{\rho}V}\] %
is a decomposition of the birational involution $\tau\in\Bir(V)$
in a sequence of elementary links (or Sarkisov links) with a
mid-point $X$ (see \cite{Co95}, \cite{CPR}, and \cite{Is96}).

Suppose that $f_{2}(x,y,z,w)$ and $f_{4}(x,y,z,w)$ are defined
over $\QQ$ and
$$
f_{3}(x,y,z,w)=\sqrt{2}g_{3}(x,y,z,w),
$$
where $g_{3}(x,y,z,t)$ is defined over $\QQ$ as well. Then the
quartic $3$-fold $V$ is defined over $\QQ(\sqrt{2})$ and not
invariant under the action of $\Gal(\QQ(\sqrt{2})\slash\QQ)$.
However, the sextic surface $S\subset\PP^3$ is given by the
equation
$$
2g^2_{3}(x,y,z,w)-4f_{2}(x,y,z,w)f_{4}(x,y,z,w)=0\subset
\PP^3=\Proj(\QQ[x,y,z,w]),
$$
which implies that the $3$-fold $X$ is defined over $\QQ$ as well.
Moreover, the $\Gal(\QQ(\sqrt{2})\slash\QQ)$-invariant part of the
group $\Cl(X)$ is $\ZZ$. Therefore, the $3$-fold $X$ is
$\QQ$-factorial and birationally super-rigid over $\QQ$ by
Theorem~A and Remark~\xref{remark-non-closed-fields}.

\end{example}

\begin{example}\label{example-on-bidegree-2-3}
Let $V$ be a smooth divisor of bidegree $(2,3)$ in
$\PP^1\times\PP^3$. The $3$-fold $V$ is then defined by the
bihomogeneous equation
\[f_3(x,y,z,w)s^2+g_3(x,y,z,w)st+h_3(x,y,z,w)t^2=0,\]
where $f_3$, $g_3$, and $h_3$ are homogeneous polynomial of degree
$3$. In addition, we denote the natural projection of $V$ to
$\PP^3$ by  $\pi:V\longrightarrow\PP^3$.   Suppose that the
polynomials $f_3$, $g_3$, and $h_3$ are general enough. The
$3$-fold $V$ then has exactly $27$ lines $C_1, C_2, \cdots,
C_{27}$ such that $-K_V\cdot C_i=0$ because the intersection
\[\{f_3(x,y,z,w)=0\}\cap\{g_3(x,y,z,w)=0\}\cap\{h_3(x,y,z,w)=0\}\] in $\PP^3$ consists of exactly $27$
points. The projection $\pi$ has degree $2$ in the outside of the
$27$ points $\pi(C_i)$. The anticanonical model
$$\Proj\left(\bigoplus_{n\geq 0} H^0\left(V, \mathcal{O}_V(-nK_V)\right)\right)$$ of $V$
is the double cover  $X$ of $\PP^3$ ramified along the nodal
sextic $S$ defined by
\[g_3^2(x,y,z,w)-4f_3(x,y,z,w)h_3(x,y,z,w)=0.\]
It has exactly $27$ simple double points each of which comes from
each line $C_i$. The morphism $\phi_{|-K_V|} : V\longrightarrow X$
given by the anticanonical system of $V$ contracts these $27$
lines to the simple double points. Therefore, it is a small
contraction and hence the double cover $X$ cannot be
$\QQ$-factorial. A generic divisor of bidegree $(2,3)$ in
$\PP^1\times\PP^3$ over $\CC$ is known to be non-rational (see
\cite{Bar84}, \cite{Ch04d}, and \cite{So02}), and hence the double
cover $X$ is also non-rational.
\end{example}

As shown in
Examples~\ref{example-on-Q-factoriality},~\ref{example-on-quartic-with-double-point},
and~\ref{example-on-bidegree-2-3}, there are non-$\QQ$-factorial
sextic double solids with $15$, $24$, and $27$ simple double
points. However, we will prove the following:

\begin{theoremB}
\label{proposition-on-Q-factoriality} Let $\pi:X\to \PP^3$ be a
double cover ramified along a  nodal sextic surface $S\subset
\PP^3$. Then the $3$-fold $X$ is $\QQ$-factorial when
$\#|\Sing(S)|\leq 14$ and it is not $\QQ$-factorial when
$\#|\Sing(S)|\geq 57$.
\end{theoremB}
Using Theorem~A with the theorem above, we immediately obtain:
\begin{corollaryB}
\label{corollary-on-non-rationality-with-small-number-of-points}
Let $\pi:X\to \PP^3$ be a double cover ramified  along a sextic
$S\subset \PP^3$ with at most $14$ simple double points. Then $X$
is birationally super-rigid. In particular, $X$ is not rational
and not birationally isomorphic to a conic bundle.
\end{corollaryB}

In \cite{CaCe82}, there are explicit constructions of sextic
surfaces in $\PP^3$ with each number of simple double points not
exceeding $64$, which give us many examples of non-rational
singular sextic double solids with at most $14$ simple double
points.

Besides the birational super-rigidity, a $\QQ$-factorial double
cover of $\PP^{3}$ ramified in a sextic nodal surface  has other
interesting properties. Implicitly the method of \cite{IsMa71} to
prove the birational (super-)rigidity also gives us information on
birational transformations to elliptic fibrations and Fano
varieties with canonical singularities.

\begin{constructionA}
\label{construction-of-elliptic-fibration} Consider a double cover
$\pi:X\to \PP^3$ ramified along a sextic $S\subset \PP^3$ with a
simple double point $o$. Let $f:W\to X$ be the blow up at the
point $o$. Then the anticanonical linear system $|-K_{W}|$ is free
and the morphism $\phi_{|-K_{W}|}:W\to \PP^{2}$ is an elliptic
fibration such that the diagram
\[\xymatrix{
&W\ar[d]_{\phi_{|-K_{W}|}}\ar@{->}[rr]^{f}&&X\ar[d]^{\pi}&\\%
&\PP^2&&\ar@{-->}[ll]_{\gamma}\PP^{3}&}\] %
is commutative, where $\gamma:\PP^{3}\dasharrow \PP^2$ is the
projection from the point $\pi(o)$.
\end{constructionA}

It is a surprise that some  double covers of $\PP^3$ ramified in
nodal sextics  can be birationally transformed into elliptic
fibrations in a way very different from the one described in
Construction~A.

\begin{constructionB}
\label{construction-of-exceptional-elliptic-fibration} Let
$\pi:X\to \PP^3$ be a double cover ramified along a sextic
$S\subset\PP^3$ such that the surface $S$ contains a line
$L\subset \PP^3$ and the line $L$ passes through exactly four
simple double points of $S$. For a general enough point $p\in X$,
there is a unique hyperplane $H_p\subset\PP^3$ containing $\pi(p)$
and $L$. The set $L\cap(C\setminus \Sing(S))$ consists of a single
point $q_p$, where $C\subset H_p$ is the quintic curve given by
$S\cap H_p=L\cup C$. The two points $\pi(p)$ and $q_p$ determine a
line $L_p$ in $\PP^3$. Define a rational map $\Xi_{L}:X\dasharrow
Grass(2,4)$ by $\Xi_{L}(p)=L_p$. The image of the map $\Xi_{L}$ is
isomorphic to $\PP^2$, hence we may assume that the map $\Xi_L$ is
a rational map of $X$ onto $\PP^2$. Obviously the map $\Xi_{L}$ is
not defined over $L$, the normalization of its general fiber is an
elliptic curve, and a resolution of indeterminacy of the map
$\Xi_{L}$ birationally transforms the $3$-fold $X$ into an
elliptic fibration.
\end{constructionB}

In this paper we will prove that these two constructions are
essentially the only ways to transform $X$ birationally into an
elliptic fibration when $X$ is $\QQ$-factorial.

\begin{theoremC}
\label{theorem-second} Let $\pi:X\to \PP^3$ be a $\QQ$-factorial
double cover ramified along a nodal sextic $S$. Suppose that we
have a birational map $\rho:X\dasharrow Y$, where $\tau:Y\to Z$ is
an elliptic fib\-ration. Then one of the following holds:
\begin{enumerate}
\item there are a simple double point $o$ on $X$ and a birational
map $\beta: \PP^2\dasharrow Z$ such that the projection $\gamma$
from the point $\pi(o)$ makes the diagram
\[ \xymatrix{
X\ar[d]_{\pi}\ar@{-->}[rr]^{\rho}&& Y\ar[d]^{\tau}\\
\PP^3\ar@{-->}[r]^{\gamma}&\PP^2\ar@{-->}[r]^{\beta}&Z}\]%
commute. \item the sextic $S$ contains a line $L\subset \PP^3$
with $\#|\Sing(S)\cap L|=4$ and there is a birational map
$\beta:\PP^2\dasharrow Z$ such that the diagram
\[ \xymatrix{
X\ar@{-->}[d]_{\Xi_{L}}\ar@{-->}[rr]^{\rho}&& Y\ar[d]^{\tau}\\
\PP^2\ar@{-->}[rr]^{\beta}&&Z}\] %
is commutative, where $\Xi_{L}$ is the rational map defined in
Construction~B.
\end{enumerate}
\end{theoremC}

In the case of one simple double point, Theorem~C was proved in
\cite{Ch01a}.

\begin{corollaryC}
\label{corollary-second} All  birational transformations of a
$\QQ$-factorial double cover of $\PP^3$ ramified along a sextic
nodal surface into elliptic fibrations\footnote{Fibrations
$\tau_1:U_1\to Z_1$ and $\tau_2:U_2\to Z_2$ can be identified if
there are birational maps $\alpha:U_1\dasharrow U_2$ and
$\beta:Z_1\dasharrow Z_2$ such that $\tau_2\circ \alpha=\beta\circ
\tau_1$ and the map $\alpha$ induces an isomorphism between
generic fibers of $\tau_1$ and $ \tau_2$.} are described by
Constructions~A and~B.
\end{corollaryC}

The following result was also proved in \cite{Ch00a}.

\begin{corollaryCC}
\label{corollary-no-elliptic-fibration} A smooth double cover $X$
of $\PP^3$ ramified along a sextic surface $S\subset \PP^3$ cannot
be birationally transformed into any elliptic fibration.
\end{corollaryCC}

\begin{remark}
\label{remark-on-double-space-with-double-line} Let $X$ be a
double cover of $\PP^{3}$ ramified in a sextic surface $S\subset
\PP^3$ such that the surface $S$ has a double line (see
\cite{Gr00}). Then the set of birational transformations of $X$
into elliptic fibrations is infinite and cannot be effectively
described (see \cite{Ch04b}).
\end{remark}

 The statement of Theorem~C is valid over an arbitrary field
$\FF$ of characteristic zero, but in Construction~A the singular
point must be defined over $\FF$  as we see in the example below.
Similarly the same has to be satisfied for Theorem~D, but the
total number of singular points on a line must be counted in
geometric sense (over the algebraic closure of $\FF$).

\begin{example}
\label{example-of-Jihun-first} Let $X$ be the double cover of
$\PP^3$ ramified in a sextic $S\subset\PP^3$ and  defined by the
equation
$$
u^{2}=x^6+xy^5+y^6+(x+y)(z^5-2zw^4)+y(z^4-2w^4)(z-3w)
$$
in $\PP(1,1,1,1,3)$. Then $X$ is smooth in the outside of $4$
simple double points given by $x=y=z^4-2w^4=0$. Hence, $X$ is
$\QQ$-factorial, birationally super-rigid, and non-rational over
$\CC$ by Theorems~A and ~B. Moreover, $x=y=0$ cuts a curve
$C\subset X$ such that $-K_{X}\cdot C=1$ and $\pi(C)\subset S$ is
a line. Therefore, $X$ can be birationally transformed over $\CC$
into exactly $5$ elliptic fibrations given by Constructions~A
and~B. However,  the $3$-fold $X$ defined over $\QQ$ is
birationally isomorphic to only one elliptic fibration given by
Construction~B.
\end{example}

Birational transformations  of other higher-dimensional algebraic
varieties into elliptic fibrations were studied in \cite{Ch00a},
\cite{Ch00b}, \cite{Ch00c}, \cite{Ch03a}, \cite{Ch03b},
\cite{Ch04a}, \cite{Ch04b}, and \cite{Ry02}. It turns out that
classification of birational transformations into elliptic
fibrations implicitly gives classification of birational
transformations into canonical Fano $3$-folds.

In the present paper we will prove the following result.

\begin{theoremD}
\label{theorem-on-canonical-Fanos} Let $\pi:X\to \PP^3$ be a
$\QQ$-factorial double cover ramified in a nodal sextic $S\subset
\PP^3$. Then $X$ is birationally isomorphic to a Fano $3$-fold
with canonical singularities that is not biregular to $X$ if and
only if the sextic $S$ contains a line $L$ passing through $5$
simple double points of the surface $S\subset \PP^3$.
\end{theoremD}

During the proof of Theorem~D, we will explicitly describe the
constructions of all possible birational transformations of sextic
double solids into Fano $3$-folds with canonical singularities.

\begin{example}
\label{example-of-Jihun-second} Let $X$ be the double cover of
$\PP^3$ ramified in a sextic $S\subset\PP^3$ and defined by the
equation
$$
u^{2}=x^6+xy^5+y^6+(x+y)(z^5-zw^4)
$$
in $\PP(1,1,1,1,3)$. Then $X$ is smooth in the outside of $5$
simple double points given by $x=y=z(z^4-w^5)=0$. For the same
reason as in Example~\xref{example-of-Jihun-first}, the double
cover $X$ is $\QQ$-factorial, birationally super-rigid, and
non-rational. As for elliptic fibrations, it can be birationally
transformed into $5$ elliptic fibrations given by Construction~A.
Also, the $3$-fold $X$ is birationally isomorphic to a unique Fano
$3$-fold with canonical singularities that is not biregular to
$X$.
\end{example}

The statements of Theorems~A, C, and~D are valid over all fields
of characteristic zero, but over fields of positive characteristic
some difficulties may occur. Indeed, the vanishing theorem of
Yu.Kawamata and E.Viehweg (see \cite{Ka82}, \cite{Vi82}) is no
longer true in positive characteristic. Even though there are some
vanishing theorems over fields of positive characteristic (see
\cite{EsVi92}, \cite{She97}), they are not applicable to our case.
A smooth resolution of indeterminacy of a birational map may fail
as well because it implicitly  uses resolution of singularities
(see \cite{Hi64}) which is completely proved  only in
characteristic zero. However, resolution of singularities for
$3$-folds is proved in \cite{Ab98} for the case of characteristic
$>5$ (see also \cite{Cos96}).

Consider the following very special example.

\begin{example}
\label{example-field-of-five-elements} Suppose that the base field
is $\FF_{5}=\ZZ\slash 5\ZZ$. Let $X$ be the double cover of
$\PP^{3}=\Proj(\FF_{5}[x,y,z,w])$ ramified along the sextic $S$
given by the equation
$$x^5y+x^4y^2+x^2y^3z-y^5z-2x^4z^2+xz^5+yz^5+x^3y^2w+2x^2y^3w-$$
$$-xyz^3w-xyz^2w^2-x^2yw^3+xy^2w^3+x^2zw^3+xyw^4+xw^5+2yw^5=0.$$
Then $X$ is smooth (see \cite{EiGr01} and \cite{GrSt}) and
$\Pic(X)\cong\ZZ$ by Lemma~3.2.2 in \cite{Do82} or Lemma~3.5 in
\cite{CPR} (see \cite{CalLy94} and \cite{Gro65}). Moreover, $X$
contains a curve $C$ given by the equations $x=y=0$ whose image in
$\PP^3$ is a line $L$ contained in the sextic $S\subset\PP^3$. For
a general enough point $p\in X$, there is a unique hyperplane
$H_p\subset\PP^3$ containing $\pi(p)$ and $L$. The residual
quintic curve $Q\subset H_p$ given by $S\cap H_p=L\cup Q$
intersects $L$ at a single point $q_p$ with
$\mult_{q_p}(Q\vert_{L})=5$.
 The two points $\pi(p)$ and $q_p$ determine a line
$L_p$ in $\PP^3$. As in Construction~B we can define a rational
map $\Psi:X\dasharrow \PP^2$ by the lines $L_p$. As we see, the
situation is almost same as that of Construction~B. We, at once,
see that a resolution of indeterminacy of the map $\Psi$
birationally transforms the $3$-fold $X$ into an elliptic
fibration.
\end{example}

Therefore, Theorem~C and even Corollary~C2 are not valid over some
fields of positive characteristic. We will however prove the
following result:

\begin{theoremE}
\label{proposition-on-finite-fields} Let $\pi:X\to \PP^3$ be a
double cover defined over a perfect field $\FF$ and ramified along
a sextic nodal surface $S\subset \PP^3$. Suppose that $X$ is
$\QQ$-factorial and $\Pic(X)\cong\ZZ$. Then $X$ is birationally
super-rigid and birational maps of $X$ into elliptic fibrations
are described by Constructions~A and~B if $\cha(\FF)>5$.
\end{theoremE}

Non-rationality and related questions like non-ruledness or
birational rigidity over fields of positive characteristic may be
interesting in the following cases:
\begin{enumerate}
\item arithmetics of algebraic varieties over
finite fields (see \cite{Es03}, \cite{LaPe00}, and \cite{Pe03});%
\item classification of varieties over fields of positive
characteristic (see \cite{Meg98} and  \cite{She97});%
\item algebro-geometric coding theory (see \cite{Bo98},
\cite{Go82}, \cite{Go88}, \cite{HeTs92}, \cite{TsVl91}, and
\cite{vdGr01});%
\item proofs of the non-rationality of certain higher-dimensional
varieties by means of reduction into fields of positive
characteristic (see \cite{ChWo04}, \cite{Ko95}, \cite{Ko96}, and
\cite{Ko00}), where even non-perfect fields may appear in some
very subtle questions as in \cite{Ko00}.
\end{enumerate}

In arithmetic geometry, it is an important and difficult problem
to measure the size of the set of rational points on a given
variety defined over a number field $\FF$. One of the most
profound works in this area is, for example, Faltings Theorem that
a smooth curve of genus at least two defined over a number field
$\FF$ has finitely many $\FF$-rational points (see \cite{Fa83}).
One of higher dimensional generalizations of the theorem is the
Weak Lang Conjecture that the set of rational points of a smooth
variety of general type defined over a number field is not Zariski
dense, which is still far away from  proofs.

A counterpart of the Weak Lang Conjecture is the conjecture that
for a smooth variety $X$ with ample $-K_X$ defined over a number
field $\FF$ there is a finite field extension of the field $\FF$
over which the set of rational points of $X$ is Zariski dense. We
can easily check that this conjecture is true for curves and
surfaces, where the condition implies that $X$ is rational over
some finite field extension. Therefore, smooth Fano $3$-folds are
the first nontrivial cases testing the conjecture.
\begin{definition}
The set of rational points of a variety $X$ defined over a number
field $\FF$ is said to be potentially dense if for some finite
field extension $\KK$ of the field $\FF$ the set of $\KK$-rational
points of $X$ is Zariski dense in $X$.
\end{definition}

Using elliptic fibrations, \cite{BoTsch98} and \cite{HaTsch00}
have proved:
\begin{theorem}
\label{theorem-on-potential-density} The set of rational points is
potentially dense on all smooth Fano $3$-folds defined over a
number field $\FF$ possibly except  double covers of $\PP^{3}$
ramified along smooth sextics.
\end{theorem}

Arithmetic properties of algebraic varieties are closely related
to their biregular and birational geometry (see \cite{BaMa90},
\cite{BaTsc95}, \cite{BaTsc96}, \cite{BaTsc98}, \cite{FrMaTsc89},
\cite{Ma66}, \cite{Ma67}, \cite{Ma72},  \cite{Ma93}, \cite{Ma95},
\cite{MaTsc93}, and \cite{MaTs86}). For example, the possible
exception appears in Theorem \xref{theorem-on-potential-density}
because smooth double covers of $\PP^{3}$ ramified in sextics are
the only smooth Fano $3$-folds that are not birationally
isomorphic to  elliptic fibrations (see \cite{IsPr99}). Besides
Fano varieties,  on several other classes of algebraic varieties
the potential density of rational points has been proved (see
\cite{BoTsch98}, \cite{BoTsch99}, and \cite{BoTsch00}).

In Section~\ref{section-potential-density} we  prove the following
result:
\begin{theoremF}
\label{theorem-on-potential-density-of-double-space} Let $\pi:X\to
\PP^3$ be a  double cover defined over a number field $\FF$ and
ramified along a sextic nodal surface $S\subset \PP^3$. If
$\Sing(X)\ne\emptyset$, then the set of rational points on $X$ are
potentially dense.
\end{theoremF}
As shown in Theorem~C, the sextic double solid can be birationally
transformed into an elliptic fibration if it has a simple double
point. Therefore, we can adopt the methods of \cite{BoTsch98} and
\cite{HaTsch00} in this case.

\vspace{5mm}

\emph{Acknowledgement.} We would like to thank V.Alexeev,
F.Bogomolov, A.Corti, M.Grinenko, V.Iskovskikh, M.Mella,
A.Pukh\-li\-kov, V.Sho\-ku\-rov, Yu.Tschinkel, and L.Wotzlaw for
helpful conversations. This work has been done during the first
author's stay at KIAS and POSTECH  in Korea. We would also like to
thank them for their hospitality. The second author has been
partially  supported by POSTECH BSRI research fund-2004.


\section{Movable log pairs and N\"other-Fano inequalities.}

To study sextic double solids we frequently use movable log pairs
introduced in \cite{Al91}.  In this section we overview their
properties and N\"other-Fano inequalities that are the most
important tools for birational (super-)rigidity.

\begin{definition}
\label{definition-of-movable-log-pair} On a variety $X$ a movable
boundary $\MMM_{X}=\sum_{i=1}^{n} a_{i}\MMM_{i}$ is a formal
finite $\QQ$-linear combination of linear systems $\MMM_i$ on $X$
such that the base locus of each $\MMM_i$ has codimension at least
two and each coefficient $a_i$ is non-negative. A movable log pair
$(X, \MMM_X)$ is a variety $X$ with a movable boundary $\MMM_X$ on
$X$.
\end{definition}

Every movable log pair can be considered as a usual log pair by
replacing each linear system   by its general element. In
particular, for a given movable log pair $(X, \MMM_{X})$ we may
handle the movable boundary $\MMM_{X}$ as an effective divisor. We
can also  consider the self-intersection $\MMM_{X}^{2}$ of
$\MMM_X$ as a well-defined effective codimension-two cycle when
$X$ is $\QQ$-factorial. We call $K_{X}+\MMM_{X}$ the log canonical
divisor of the movable log pair $(X, \MMM_{X})$. Throughout the
rest of this section, we will assume that log canonical divisors
are $\QQ$-Cartier divisors.

\begin{definition}
\label{birationally-equivalent-movable-log-pairs} Movable log
pairs $(X, \MMM_{X})$ and $(Y, \MMM_{Y})$ are birationally
equivalent if there is a birational map $\rho:X\dasharrow Y$ such
that $\MMM_{Y}=\rho(\MMM_{X})$.
\end{definition}

The notions such as discrepancies,  (log) terminality, and (log)
cano\-ni\-city can be defined for movable log pairs as for usual
log pairs (see \cite{KMM}).

\begin{definition}
\label{definition-of-canonical-singularities} A movable log pair
$(X, \MMM_{X})$ has canonical (terminal, resp.) singularities if
for every birational morphism $f:W\to X$ each discrepancy $a(X,
\MMM_{X}, E)$ in
$$
K_{W}+f^{-1}(\MMM_{X})\sim_{\QQ} f^{*}(K_{X}+\MMM_{X})+\sum_{E
\mbox{{\scriptsize : $f$-exceptional divisor}}}a(X, \MMM_{X}, E)E
$$
is non-negative (positive, resp.).
\end{definition}

\begin{example}
\label{example-of-movable-log-pairs-with-reduced-boundary} Let
$\MMM$ be a linear system on a 3-fold $X$ with no fixed
components. Then the log pair $(X, \MMM)$ has terminal
singularities if and only if the linear system $\MMM$ has only
isolated simple base points which are smooth points on the 3-fold
$X$.
\end{example}

 Log Minimal Model Program  holds good
for  three-dimensional movable log pairs with canonical (terminal)
singularities (see \cite{Al91} and \cite{KMM}). In particular, it
preserves their canonicity (terminality).

Every movable log pair is birationally equivalent to a movable log
pair with canonical or terminal singularities. Away from the base
loci of the components of its boundary, the singularities of a
movable log pair coincide with those of its variety.

\begin{definition}
\label{center-of-canonical-singularities} A proper irreducible
subvariety $Y\subset X$ is called a center of the canonical
singularities of a movable log pair $(X, \MMM_{X})$ if there are a
birational morphism $f:W\to X$ and an $f$-ex\-cep\-tional divisor
$E\subset W$ such that the discrepancy $a(X, \MMM_{X}, E)\leq 0$
and $f(E)=Y$. The set of all the centers of the canonical
sin\-gu\-la\-ri\-ties of the movable log pair $(X, \MMM_{X})$ will
be denoted by $\CCSS(X, \MMM_{X})$. 
\end{definition}


Note that  a log pair $(X, \MMM_{X})$ is terminal if and only if
$\CCSS(X, \MMM_{X})=\emptyset$.

Let $(X, \MMM_{X})$ be a movable log pair and $Z\subset X$ be a
proper irreducible subvariety such that $X$ is smooth along the
subvariety $Z$. Then elementary properties of blow ups along
smooth subvarieties of smooth varieties imply that
$$
Z\in \CCSS(X, \MMM_{X})\Rightarrow \mult_{Z}(\MMM_{X})\geq 1
$$
and in the case when $\codim(Z\subset X)=2$ we have
$$
Z\in \CCSS(X, \MMM_{X})\iff \mult_{Z}(\MMM_{X})\geq 1.
$$

For a movable log pair $(X,\MMM_X)$ we consider a birational
morphism $f:W\to X$ such that the log pair $(W,
\MMM_W:=f^{-1}(\MMM_X))$ has canonical singularities.
\begin{definition}
\label{Kodaira-dimension} The number $\kappa(X,
\MMM_{X})=\dim(\phi_{|nm(K_{W}+\MMM_{W})|}(W))$ for $n\gg 0$ is
called the Ko\-dai\-ra dimension of the movable log pair $(X,
\MMM_{X})$, where $m$ is a natural number such that
$m(K_{W}+\MMM_{W})$ is a Cartier divisor.
 When
$|nm(K_{W}+\MMM_{W})|=\emptyset$ for all $n\in \NN$, the Kodaira
dimension  $\kappa(X, \MMM_{X})$ is defined to be $-\infty$.
\end{definition}

\begin{proposition}
\label{independence-of-Kodaira-dimension} The Kodaira dimension of
a movable log pair is well-defined. In particular, it does not
depend on the choice of the birationally equivalent movable log
pair with canonical singularities.
\end{proposition}
\begin{proof}
Let $(X, \MMM_{X})$ and $(Y, \MMM_{Y})$ be movable log pairs with
canonical singularities such that there is a birational map
$\rho:Y\dasharrow X$ with $\MMM_{X}=\rho(\MMM_{Y})$. Choose
positive integers $m_{X}$ and $m_{Y}$ such that both
$m_X(K_{X}+\MMM_{X})$ and $m_Y(K_{Y}+\MMM_{Y})$ are Cartier
divisors. We must show that either
\[|nm_{X}(K_{X}+\MMM_{X})|=|nm_{Y}(K_{Y}+\MMM_{Y})|=\emptyset \mbox{ for all } n\in\NN \] or
\[\dim(\phi_{|nm_{X}(K_{X}+\MMM_{X})|}(X))=\dim(\phi_{|nm_{Y}(K_{Y}+\MMM_{Y})|}(Y))
\mbox{ for }n\gg 0.\] We consider a Hironaka hut of
$\rho:Y\dasharrow X$, \emph{i.e.}, a smooth variety $W$ with
birational morphisms $g:W\to X$ and $f:W\to Y$ such that the
diagram
\[\xymatrix{
&W\ar@{->}[dr]^{g}\ar@{->}[dl]_{f}&\\%
Y\ar@{-->}[rr]^{\rho}&&X}\] %
commutes.  We then obtain
$$
K_{W}+\MMM_{W}\sim_{\QQ}
g^{*}(K_{X}+\MMM_{X})+\Sigma_{X}\sim_{\QQ}
f^{*}(K_{Y}+\MMM_{Y})+\Sigma_{Y},
$$
where $\MMM_{W}=g^{-1}(\MMM_{X})$, $\Sigma_{X}$ and $\Sigma_{Y}$
are the exceptional divisors of $g$ and $f$ respectively. Because
the movable log pairs $(X, \MMM_{X})$ and $(Y, \MMM_{Y})$ have
canonical singularities, the exceptional divisors $\Sigma_{X}$ and
$\Sigma_{Y}$ are effective and hence the linear systems
$|n(K_{W}+\MMM_{W})|$, $|g^{*}(n(K_{X}+\MMM_{X}))|$, and
$|f^{*}(n(K_{Y}+\MMM_{Y}))|$ have the same dimension for a big and
divisible enough natural number $n$. Moreover, if these linear
systems are not empty, then we have
$$
\phi_{|n(K_{W}+\MMM_{W})|}=\phi_{|g^{*}(n(K_{X}+\MMM_{X}))|}=\phi_{|f^{*}(n(K_{Y}+\MMM_{Y}))|},
$$
which implies the claim.
\end{proof}

By definition, the Kodaira dimension of a movable log pair is a
birational invariant and a non-decreasing function of the
coefficients of the movable boundary.

\begin{definition}
\label{canonical-model} A movable log pair $(V, \MMM_{V})$ is
called
 a canonical model of a movable log pair $(X, \MMM_{X})$ if
there is a birational map $\psi:X\dasharrow V$ such that
$\MMM_{V}=\psi(\MMM_{X})$, the movable log pair $(V, \MMM_{V})$
has canonical singularities, and the divisor $K_{V}+\MMM_{V}$ is
ample.
\end{definition}

\begin{proposition}
\label{uniqueness-of-canonical-model} A canonical model of a
movable log pair is unique if it exists.
\end{proposition}
\begin{proof}
Let $(X, \MMM_{X})$ and $(Y, \MMM_{Y})$ be canonical models such
that there is a birational map  $\rho:Y\dasharrow X$ with
$\MMM_{X}=\rho(\MMM_{Y})$. Take  a smooth variety $W$ with
birational morphisms $g:W\to X$ and $f:W\to Y$ such that the
diagram
\[\xymatrix{
&W\ar@{->}[dr]^{g}\ar@{->}[dl]_{f}&\\%
Y\ar@{-->}[rr]^{\rho}&&X}\] %
commutes. We have
$$
K_{W}+\MMM_{W}\sim_{\QQ}
g^{*}(K_{X}+\MMM_{X})+\Sigma_{X}\sim_{\QQ}
f^{*}(K_{Y}+\MMM_{Y})+\Sigma_{Y},
$$
where $\MMM_{W}=g^{-1}(\MMM_{X})=f^{-1}(\MMM_{Y})$, $\Sigma_{X}$
and $\Sigma_{Y}$ are the exceptional divisors of birational
morphisms $g$ and $f$ respectively. Let $n\in \NN$ be a big and
divisible enough number such that $n(K_{W}+\MMM_{W})$,
$n(K_{X}+\MMM_{X})$, and $n(K_{Y}+\MMM_{Y})$ are Cartier divisors.
For the same reason as in the proof of
Proposition~\xref{independence-of-Kodaira-dimension} we obtain
$$
\phi_{|n(K_{W}+\MMM_{W})|}=\phi_{|g^{*}(n(K_{X}+\MMM_{X}))|}=\phi_{|f^{*}(n(K_{Y}+\MMM_{Y}))|}
$$
Therefore, the birational map  $\rho$ is an isomorphism because
$K_{X}+\MMM_{X}$ and $K_{Y}+\MMM_{Y}$ are ample.
\end{proof}

The existence of the canonical model of a movable log pair implies
that its Kodaira dimension equals to the dimension of the variety.

N\"other-Fano inequalities can be immediately reinterpreted in
terms of canonical singularities of movable log pairs. For
reader's understanding, we give the theorems and their proofs on
the relation between singularities of movable log pairs and
birational (super-)rigidity. In addition, with del Pezzo surfaces
of Picard number $1$ defined over non-closed fields, we
demonstrate how to apply the theorems, which is so simple that one
can easily understand.

The following result is known as a classical N\"other-Fano
inequality.
\begin{theorem}
\label{theorem-Nother-Fano-inequality} Let $X$ be a terminal
$\QQ$-factorial Fano variety with $\Pic(X)\cong \ZZ$. If every
movable log pair $(X, \MMM_{X})$ with $K_{X}+\MMM_{X}\sim_\QQ 0$
has canonical singularities, then $X$ is birationally super-rigid.
\end{theorem}

\begin{proof}
Suppose that there is a birational map $\rho:X\dasharrow V$ such
that $V$ is a Fano variety with $\QQ$-factorial terminal
singularities and $\Pic(V)\cong \ZZ$. We are to show that $\rho$
is an isomorphism. Let $\MMM_{V}=r |-nK_{V}|$ and
$\MMM_{X}=\rho^{-1}(\MMM_{V})$ for a natural number $n\gg 0$ and a
rational number $r>0$ such that $K_{X}+\MMM_{X}\sim_{\QQ} 0$.
Because $|-nK_V|$ is free for $n\gg 0$ and $V$ has at worst
terminal singularities, the log pair $(V, \MMM_V)$ has terminal
singularities. In addition, the equality
$$\kappa(X, \MMM_{X})=\kappa(V, \MMM_{V})=0.$$
implies that the divisor $K_V+\MMM_V$ is nef; otherwise  the
Kodaira dimension $\kappa(V, \MMM_{V})$ would be $-\infty$.

Let $f:W\to X$ be a birational morphism of a smooth variety $W$
such that $g=\rho\circ f$ is a morphism. Then
\begin{align*}
K_{W}+\MMM_{W}&= f^{*}(K_{X}+\MMM_{X})+\sum_{i=1}^{l_1}a(X,
\MMM_{X}, F_{i})F_{i}+\sum_{k=1}^{m}a(X, \MMM_{X}, E_{k})E_{k}\\
& =g^{*}(K_{V}+\MMM_{V})+\sum_{j=1}^{l_2}a(V, \MMM_{V},
G_{j})G_{j}+\sum_{k=1}^{m}a(V, \MMM_{V}, E_{k})E_{k},
\end{align*}
where $\MMM_{W}=f^{-1}(\MMM_{X})$, each divisor $F_{i}$ is
$f$-exceptional but not $g$-exceptional, each divisor $G_{j}$ is
$g$-exceptional but not $f$-exceptional, and each $E_k$ is both
$f$-exceptional and $g$-exceptional. Applying Lemma 2.19 in
\cite{Ko91}, we obtain
$$
a(X, \MMM_{X}, E_k)=a(V, \MMM_{V}, E_k)
$$
for each $k$ and we see that there is no $g$-exceptional but not
$f$-exceptional divisor, \emph{i.e.}, $l_2=0$ because the log pair
$(V, \MMM_V)$ has terminal singularities.  Furthermore, there
exits no $f$-exceptional but not $g$-exceptional divisor,
\emph{i.e.}, $l_1=0$ because the  Picard numbers of $V$ and $X$
are the same. Therefore, 
the log
pair $(X, \MMM_{X})$ has at worst terminal singularities. For some
$d\in \QQ_{>1}$, both the movable log pairs $(X, d \MMM_{X})$ and
$(V, d\MMM_{V})$ are canonical models. Hence, $\rho$ is an
isomorphism by Proposition~\xref{uniqueness-of-canonical-model}.

We now suppose that we have a birational map $\chi:X\dasharrow Y$
of $X$ into a fibration $\tau:Y\to Z$, where $Y$ is smooth and a
general fiber of $\tau$ is a smooth variety of Kodaira dimension
$-\infty$. Let $\MMM_{Y}=c|\tau^{*}(H)|$ and
$\MMM_{X}=\chi^{-1}(\MMM_{Y})$, where $H$ is a very ample divisor
on $Z$ and  $c$ is a positive rational number such that
$K_{X}+\MMM_{X}\sim_{\QQ} 0$. Then the Kodaira dimension
$\kappa(X, \MMM_{X})$ is zero because the log pair $(X, \MMM_{X})$
has at worst canonical singularities and $K_{X}+\MMM_{X}\sim_{\QQ}
0$. However, the Kodaira dimension $\kappa(Y, \MMM_{Y})=-\infty.$
This contradiction completes the proof.
\end{proof}

The proof of Theorem~\xref{theorem-Nother-Fano-inequality} shows a
condition for the Fano variety $X$ to be birationally rigid as
follows:

\begin{corollary}
\label{corollary-on-birational-rigidity} Let $X$ be a terminal
$\QQ$-factorial Fano variety with $\Pic(X)\cong \ZZ$. Suppose that
for every movable log pair $(X, \MMM_{X})$ with
$K_{X}+\MMM_{X}\sim_\QQ 0$ either the singularities of the log
pair $(X, \MMM_{X})$ are canonical or the divisor
$-(K_{X}+\rho(\MMM_{X}))$ is ample for some birational
automorphism $\rho\in\Bir(X)$. Then $X$ is birationally rigid.
\end{corollary}

Log Minimal Model Program tells us that the condition in
Theorem~\xref{theorem-Nother-Fano-inequality} is a necessary and
sufficient one for $X$ to be birationally super-rigid.

\begin{proposition}
\label{proposition-Nother-Fano-inequality-for-threefolds} Let $X$
be a terminal $\QQ$-factorial Fano $3$-fold with $\Pic(X)\cong
\ZZ$. Then $X$ is birationally super-rigid if and only if every
movable log pair $(X, \MMM_{X})$ with $K_{X}+\MMM_{X}\sim_\QQ 0$
has at worst canonical singularities.
\end{proposition}

\begin{proof}
Suppose that $X$ is birationally super-rigid. In addition, we
suppose that there is a movable log pair $(X, \MMM_{X})$ with
non-canonical singularities such that $K_{X}+\MMM_{X}\sim_\QQ 0$.
Let $f:W\to X$ be a birational morphism such that the log pair
$(W, \MMM_{W}:=f^{-1}(\MMM_{X}))$ has canonical singularities.
Then
$$
K_{W}+\MMM_{W}=
f^{*}(K_{X}+\MMM_{X})+\sum_{i=1}^{k}a(X, \MMM_{X}, E_{i})E_{i}\sim_\QQ %
\sum_{i=1}^{k}a(X, \MMM_{X}, E_{i})E_{i},
$$
where $E_{i}$ is an $f$-exceptional divisor and $a(X,\MMM_{X},
E_{j})<0$ for some $j$.

Applying relative Log Minimal Model Program to the log pair $(W,
\MMM_{W})$ over $X$ we may assume $K_{W}+\MMM_{W}$ is $f$-nef.
Then, Lemma 2.19 in \cite{Ko91} immediately implies that $a(X,
\MMM_{X}, E_{i})\le 0$ for all $i$. Log Minimal Model Program for
$(W, \MMM_{W})$ gives a birational map $\rho$ of $W$ into  a Mori
fibration space $Y$, \emph{i.e.}, a fibration  $\pi:Y\to Z$ such
that $-K_{Y}$ is $\pi$-ample, the variety $Y$ has $\QQ$-factorial
terminal singularities, and $\Pic(Y/Z)\cong \ZZ$. However, the
birational map $\rho\circ f^{-1}$ is not an isomorphism.
\end{proof}

Despite its formal appearance,
Theorem~\xref{theorem-Nother-Fano-inequality} can be effectively
applied in many different cases. For example, the following result
in \cite{Ma66} and \cite{Ma67} is an application of
Theorem~\xref{theorem-Nother-Fano-inequality}.

\begin{theorem}
\label{theorem-of-Manin} Let $X$ be a smooth del Pezzo surface
 defined over a
perfect field $\FF$ with $\Pic(X)\cong \ZZ$ and $K_{X}^{2}\le 3$.
Then $X$ is birationally rigid and non-rational over $\FF$.
\end{theorem}

\begin{proof}
Suppose that $X$ is not birationally rigid. Then there is a
movable log pair $(X, \MMM_{X})$ defined over $\FF$ such that
$K_{X}+\MMM_{X}\sim_\QQ 0$ and that is not canonical at some
smooth point $o\in X$. Therefore, $\mult_{o}(\MMM_{X})>1$ and
$$
3\geq K_{X}^{2}=M^{2}_{X}\geq mult^{2}_{o}(\MMM_{X})\deg(o\otimes
{\bar \FF}) >\deg(o\otimes {\bar \FF}),
$$
where ${\bar \FF}$ is the algebraic closure of the field $\FF$. In
 the case $K_{X}^{2}=1$, the strict inequality is a contradiction. Moreover, if
$K_{X}^{2}=2$, then the point $o$ is defined over $\FF$ and if
$K_{X}^{3}=3$, then the point $o$ splits into no more than  two
points over the field ${\bar \FF}$.

Suppose that $K_X^2$ is either $2$ or $3$. Let $f:V\to X$ be the
blow up at the point $o$. Then
$$
K_{V}^{2}=K_{X}^{2}-\deg(o\otimes {\bar \FF})
$$
and $V$ is a smooth del Pezzo surface because $\Pic(X)=\ZZ$, the
inequality $\mult_{o}(\MMM_{X})>1$ holds, and the boundary
$\MMM_{X}$ is movable. The double cover $\phi_{|-K_{V}|}$ induces
an involution $\tau\in \Bir(X)$ that is classically known as
Bertini or Geizer involution. Simple calculations show the
ampleness of divisor $-(K_{X}+\tau(\MMM_{X}))$, which contradicts
Corollary~\xref{corollary-on-birational-rigidity}.
\end{proof}

The proofs of Theorems~\xref{theorem-Nother-Fano-inequality} and
\xref{theorem-of-Manin} and Lemma 5.3.1 in \cite{Ko00} imply that
a result similar to Theorem~\xref{theorem-of-Manin} holds over a
non-perfect field as well. Indeed, one can prove that a
non-singular del Pezzo surface $X$  defined over non-perfect field
$\FF$ is non-rational over $\FF$ and is not birationally
isomorphic over $\FF$ to any non-singular del Pezzo surface $Y$
with $\Pic(Y)=\ZZ$, which is smooth in codimension one, if
$\Pic(X)\cong \ZZ$ and $K_{X}^{2}\le 3$.

Most applications of
Theorems~\xref{theorem-Nother-Fano-inequality} have the pattern of
the proof of Theorem~\xref{theorem-of-Manin} implicitly.

The following result can be considered as a weak N\"other-Fano
inequality.

\begin{theorem}
\label{theorem-on-elliptic-fibrations} Let $X$ be a terminal
$\QQ$-factorial Fano variety with $\Pic(X)\cong \ZZ$,
$\rho:X\dasharrow Y$  a birational map, and $\pi:Y\to Z$ a
fibration. Suppose that a general enough fiber of $\pi$ is a
smooth variety of Kodaira dimension zero. Then the singularities
of the movable log pair $(X, \MMM_{X})$ is not terminal,
where $\MMM_{X}=r\rho^{-1}(|\pi^{*}(H)|)$ for a very ample divisor
$H$ on $Z$ and $r\in\QQ_{>0}$ such that $K_{X}+\MMM_{X}\sim_\QQ
0$.
\end{theorem}

\begin{proof}
Suppose $\CCSS(X, \MMM_{X})=\emptyset$. %
Let $\MMM_{Y}=r|\pi^{*}(H)|$. Then we see
$$
\kappa(X, c \MMM_{X})=\kappa(Y, c \MMM_{Y})\le\dim(Z)<\dim(X).
$$
However, $\CCSS(X, c \MMM_{X})=\emptyset$ for small $c>1$ and
hence $\kappa(X, c \MMM_{X})=\dim(X)$, which is a contradiction.
\end{proof}

The easy result below shows how to apply
Theorem~\xref{theorem-on-elliptic-fibrations}.

\begin{proposition}
\label{proposition-on-elliptic-fibrations} Let $X$ be a smooth del
Pezzo surface of degree one with $\Pic(X)\cong \ZZ$ defined over a
perfect field $\FF$ and $o$  the unique base point of the
anticanonical linear system of the surface $X$. Let
$\rho:X\dasharrow Y$ be a birational map, where $Y$ is a smooth
surface. Suppose that $\pi:Y\to Z$ is a relatively minimal
elliptic fibration with connected fibers such that a general
enough fiber of  $\pi$ is smooth. Then the birational map $\rho$
is the  blow up at some $\FF$-rational point $p$ on the del Pezzo
surface $X$ and the morphism $\pi$ is induced by $|-nK_Y|$ for
some $n\in \NN$. Furthermore, $p\in {\hat C}$ and the equality
$p^{n}=\id_{{\hat C}}$ holds, where ${\hat C}$ is the smooth part
of the unique curve $C$ of arithmetic genus one in $|-K_{X}|$
passing though the point $p$ and considered as a group scheme with
the identity $\id_{{\hat C}}=o$.
\end{proposition}

\begin{proof}
Let $\MMM_{X}=c\rho^{-1}(|\pi^{*}(H)|)$, where $H$ is a very ample
on curve $Z$ and $c\in \QQ_{>0}$, such that the equivalence
$K_{X}+\MMM_{X}\sim_\QQ 0$ holds. Then the set $\CCSS(X,
\MMM_{X})$ contains a point $p$ on the surface $X$ by
Theorem~\xref{theorem-on-elliptic-fibrations}. In particular,
$\mult_{p}(\MMM_{X})\geq 1$, but
$$
1=K_{X}^{2}=M^{2}_{X}\geq mult^{2}_{p}(\MMM_{X})\deg(p\otimes
{\bar \FF})\geq \deg(p\otimes {\bar \FF})\geq 1,
$$
where ${\bar \FF}$ is the algebraic closure of the field $\FF$.
Hence, $\mult_{p}(\MMM_{X})=1$ and the point $p$ is defined over
the field $\FF$. Let $f:V\to X$ be the  blow up at the point $p$.
Then $K_{V}^{2}=0$ and
$$
-K_{V}\sim_\QQ \MMM_{V}=f^{-1}(\MMM_{X}),
$$
which implies that the linear system $|-rK_{V}|$ is free for a
natural number $r\gg 0$. The morphism $\phi_{|-rK_{V}|}$ is a
relatively minimal elliptic fibration and $\MMM_{V}\cdot E=0$ for
a general enough fiber $E$ of the elliptic fibration
$\phi_{|-rK_{V}|}$. Therefore the linear system $(\rho\circ
f)^{-1}(|\pi^{*}(H)|)$ is contained in the fibers of the fibration
$\phi_{|-rK_{V}|}$. Relative minimality of the fibrations $\pi$
and $\phi_{|-rK_{V}|}$ implies $\rho\circ f$ is an isomorphism.

Suppose $p\ne o$. Let $C\in |-K_{X}|$ be a curve passing through
$p$. Because
$$
1=K_{X}^{2}=C\cdot \MMM_{X}\geq %
\mult_{p}(\MMM_{X})\mult_{p}(C)=\mult_{p}(C)\geq 1,
$$
the curve $C$ is smooth at the point $p$. Let ${\tilde
C}=f^{-1}(C)\sim -K_{V}$. Then $h^{0}(V,\OOO_V({\tilde C}))=1$ and
the  curve ${\tilde C}$ is $\Gal({\bar \FF}\slash \FF)$-invariant.
In particular, the curve $Z$ has a $\FF$-point
$\phi_{|-rK_{V}|}({\tilde C})$ and we have $Z\cong \PP^{1}$. Take
the smallest natural $n$ such that $h^{0}(V,\OOO_V({n\tilde
C}))>1$. The exact sequence
$$
0\to \OOO_{V}((n-1){\tilde C})\to \OOO_{V}(n{\tilde C})\to
\OOO_{{\tilde C}}(n{\tilde C}\vert_{\tilde C})\to 0
$$
implies $h^{0}(C, \OOO_{C}(n(p-o)))=h^{0}(\tilde{C},\OOO_{{\tilde
C}}(n{\tilde C}\vert_{\tilde C}))\ne 0$, which implies the claim.
\end{proof}

\begin{corollary}
\label{corollary-on-elliptic-fibrations} Let $X$ be a terminal
$\QQ$-factorial Fano variety with $\Pic(X)\cong \ZZ$  such that
every movable log pair $(X, \MMM_{X})$ with
$K_{X}+\MMM_{X}\sim_\QQ 0$ has terminal singularities. Then $X$ is
not birationally isomorphic to a fibration of varieties of Kodaira
dimension zero.
\end{corollary}

Unfortunately Corollary~\xref{corollary-on-elliptic-fibrations} is
almost impossible to use. As far as we know, there are no known
examples of Fano varieties that are not birationally isomorphic to
fibrations of varieties of Kodaira dimension zero. The only known
example of a rationally connected variety that can not be
birationally transformed into a fibration of varieties of Kodaira
dimension zero is a conic bundle with a big enough discriminant
locus in \cite{Ch04a}.

\begin{theorem}
\label{theorem-on-canonical-Fano-varieties} Let $X$ be a terminal
$\QQ$-factorial Fano variety with $\Pic(X)\cong \ZZ$ and
$\rho:X\dasharrow Y$ be a non-biregular birational map onto a Fano
variety $Y$ with canonical singularities. Then
$K_{X}+\MMM_{X}\sim_\QQ 0$ and
$$
\CCSS(X, \MMM_{X})\ne\emptyset,
$$
where $\MMM_{X}={\frac {1} {n}}\rho^{-1}(|-nK_{Y}|)$ for a natural
number $n\gg 0$.
\end{theorem}

\begin{proof}
Let $\MMM_{Y}={\frac {1} {n}}|-nK_{Y}|$. We then see
$$
\kappa(X, \MMM_{X})=\kappa(Y, \MMM_{Y})=0,
$$
which implies $K_{X}+\MMM_{X}\sim_\QQ 0$. Suppose $\CCSS(X,
\MMM_{X})=\emptyset$.  Both the log pair $(X, r \MMM_{X})$ and
$(Y, r \MMM_{Y})$ are canonical models for a rational number $r>1$
sufficiently close to $1$. It is a contradiction that $\rho$ is an
isomorphism by Proposition~\xref{uniqueness-of-canonical-model}.
\end{proof}

The following easy result shows how to apply
Theorem~\xref{theorem-on-canonical-Fano-varieties}.

\begin{proposition}
\label{proposition-on-del-Pezzo-surface-of-degree-one} Let $X$ be
a smooth del Pezzo surface of degree one with $\Pic(X)\cong \ZZ$
defined over an arbitrary perfect field $\FF$. Then the surface
$X$ is not birationally isomorphic to a del Pezzo surface with du
Val singularities which is not isomorphic to the surface $X$.
\end{proposition}

\begin{proof}
Let $\rho:X\dasharrow Y$ be a birational map over the field $\FF$
and $\MMM_{X}={\frac {1} {n}}\rho^{-1}(|-nK_{Y}|)$ for a natural
number $n\gg 0$, where $Y$ is a del Pezzo surface with du Val
singularities and $\rho$ is not an isomorphism. Then
$K_{X}+\MMM_{X}\sim_\QQ 0$ and $\CCSS(X, \MMM_{X})$ contains some
smooth point $o$ on the del Pezzo surface $X$ by
Theorem~\xref{theorem-on-canonical-Fano-varieties}. In particular,
$\mult_{o}(\MMM_{X})\geq 1$, but
$$
1=K_{X}^{2}=M^{2}_{X}\geq mult^{2}_{o}(\MMM_{X})\deg(o\otimes
{\bar \FF})\geq \deg(o\otimes {\bar \FF})\geq 1,
$$
where ${\bar \FF}$ is the algebraic closure of the field $\FF$.
Hence, $\mult_{o}(\MMM_{X})=1$ and the point $o$ is defined over
the field $\FF$. Let $f:V\to X$ be the blow up at the point $o$.
Then $K_{V}^{2}=0$ and
$$
-K_{V}\sim_\QQ \MMM_{V}=f^{-1}(\MMM_{X}),
$$
which implies freeness of the linear system $|-rK_{V}|$ for a
natural number $r\gg 0$. The morphism $\phi_{|-rK_{V}|}$ is an
elliptic fibration and $\MMM_{V}\cdot E=0$ for a general enough
fiber $E$ of $\phi_{|-rK_{V}|}$. Therefore, the linear system
$(\rho\circ f)^{-1}(|-nK_{Y}|)$ is compounded from a pencil, which
is impossible.
\end{proof}

The paper \cite{IsMa71} of V.~Iskovskikh and Yu.~Manin  was based
on the idea of G. Fano that can be summarized by N\"other-Fano
inequalities. Since 1971 the method of V.~Iskovskikh and Yu.~Manin
has evolved to show birational rigidity of various Fano varieties.
Recently, Shokurov's connectedness principle improved the method
so that one can extremely simplify the proof of the result of
V.~Iskovskikh and Yu.~Manin (see \cite{Co00}). Furthermore, it
also made it possible to prove the birational super-rigidity of
smooth hypersurfaces of degree $n$ in $\PP^n$, $n\geq 4$ (see
\cite{Pu02}). In what follows we will explain Shokurov's
connectedness principle and how it can be applied to birational
rigidity.

 Movable boundaries always can be
considered as effective divisors and movable log pairs as usual
log pairs. Therefore, we may use compound log pairs that contains
both movable and fixed components. From now, we will not assume
any restrictions on the coefficients of boundaries. In particular,
boundaries may not be effective unless otherwise stated.

\begin{definition}
\label{log-pull-back} A log pair $(V, B^{V})$ is called the log
pull back of a log pair $(X, B_{X})$ with respect to a birational
morphism $f:V\to X$ if
$$
B^{V}=f^{-1}(B_{X})-\sum_{i=1}^{n}a(X, B_{X}, E_{i})E_{i},
$$
where $a(X, B_{X}, E_{i})$ is  the discrepancy of an
$f$-exceptional divisor $E_i$ over $(X, B_X)$. In particular, it
satisfies $K_{V}+B^{V}\sim_\QQ f^{*}(K_{X}+B_{X})$
\end{definition}

\begin{definition}
\label{center-of-log-canonical-singularities} A proper irreducible
subvariety $Y\subset X$ is  called a center of the log canonical
singularities of $(X, B_{X})$ if there are a birational morphism
$f:W\to X$ and a divisor $E\subset W$ such that $E$ is contained
in the support of the effective part of the divisor $\lfloor
B^{W}\rfloor$ and $f(E)=Y$. The set of all the centers of the log
canonical sin\-gu\-la\-ri\-ties of a log pair $(X, B_{X})$ will be
denoted by $\LLCCSS(X, B_{X})$. In addition, the union of all the
centers of log canonical singularities of $(X, \MMM_{X})$ will be
denoted by
 $LCS(X, B_{X})$.
\end{definition}

Consider a log pair $(X, B_{X})$, where $B_X=\sum_{i=1}^ka_iB_i$
is effective and $B_i$'s are prime divisors on $X$. Choose a
birational morphism $f:Y\to X$ such that $Y$ is smooth and the
union of all the proper transforms of the divisors $B_{i}$  and
all $f$-exceptional divisors forms a divisor with simple normal
crossing. The morphism $f$ is called a log resolution of the log
pair $(X, B_{X})$. By definition, the equality
$$
K_{Y}+B^{Y}\sim_\QQ f^{*}(K_{X}+B_{X})
$$
holds, where  $(Y, B^{Y})$ is the log pull back of  the log pair
$(X, B_{X})$ with respect to the birational morphism $f$.

\begin{definition}
\label{log-canonical-singularities-subscheme} The subscheme
$\LLL(X, B_{X})$ associated with the ideal sheaf $\III(X,
B_{X})=f_{*}(\mathcal{O}_Y(\lceil -B^{Y}\rceil))$ is called  the
log canonical singularity subscheme of the log pair $(X, B_{X})$.
\end{definition}

The support of the  subscheme $\LLL(X, B_{X})$ is exactly the
locus of $LCS(X, B_{X})$. The following result is called Shokurov
vanishing (see \cite{Sh92}).

\begin{theorem}
\label{vanishing-theorem-of-Shokurov} Let $(X, B_{X})$ be a log
pair with an effective divisor $B_{X}$. Suppose that there is a
nef and big $\QQ$-divisor $H$ on $X$ such that $D=K_{X}+B_{X}+H$
is Cartier. Then $H^{i}(X, \III(X, B_{X})\otimes \OOO_X(D))=0$ for
$i>0$.
\end{theorem}

\begin{proof}
Let $f:W\longrightarrow X$ be a log resolution of $(X,B_X)$.
Because $f^*H$ is nef and big on $W$ and $f^*D=K_W+B^W+f^*H$, we
obtain
$$R^{i}f_{*}(f^*\OOO_X(D)\otimes\OOO_W(\lceil -B^{W}\rceil))=0$$
for $i>0$  from relative Kawamata-Viehweg vanishing (see
\cite{Ka82} and \cite{Vi82}). The degeneration of local-to-global
spectral sequence and
$$R^{0}f_{*}(f^{*}\OOO_X(D)\otimes\OOO_W(\lceil -B^{W}\rceil))=\III(X, B_{X})\otimes \OOO_X(D)$$
imply that for all $i$
$$H^{i}(X, \III (X, B_{X})\otimes \OOO_X(D))=H^{i}(W, f^*\OOO_X(D)\otimes\OOO_W(\lceil -B^{W}\rceil)),$$
while $H^{i}(W, f^*\OOO_X(D)\otimes\OOO_W(\lceil -B^{W}\rceil))=0$
for $i>0$ by Kawamata-Viehweg vanishing.
\end{proof}

Consider the following application of
Theorem~\xref{vanishing-theorem-of-Shokurov}.

\begin{lemma}
\label{lemma-on-quadric}  Let $V$ be a variety isomorphic to
$\PP^{1}\times \PP^{1}$. Let $B_{V}$ be an effective $\QQ$-divisor
on $V$ of type $(a, b)$, where $a$ and $b\in\QQ\cap [0, 1)$. Then
$\LLCCSS(V, B_{V})=\emptyset$.
\end{lemma}

\begin{proof}
Intersecting the boundary $B_{V}$ with the rulings of $V$, we see
that the set $\LLCCSS(V, B_{V})$ does not contain   a curve on
$V$. Suppose that the set $\LLCCSS(V, B_{V})$ contains a point
$o$.  There is a $\QQ$-divisor $H$ on $V$ of type $(1-a, 1-b)$
such that the divisor
$$D=K_{V}+B_{V}+H$$
is Cartier. Since the divisor $H$ is ample,
Theorem~\xref{vanishing-theorem-of-Shokurov} implies the sequence
$$H^{0}(V,\OOO_{V}(D))\to H^{0}(\LLL(V, B_V),\OOO_{\LLL(V, B_{V})}(D))\to 0$$
is exact. However,  $H^{0}(V, \OOO_{V}(D))=0$, which  is a
contradiction.
\end{proof}

For every Cartier divisor $D$ on $X$, the sequence
$$
0\to \III(X, B_{X})\otimes D\to \OOO_{X}(D)\to \OOO_{\LLL(X,
B_{X})}(D)\to 0
$$
is exact and Theorem~\xref{vanishing-theorem-of-Shokurov} implies
the following two connectedness theorems of V.Shokurov.

\begin{theorem}
\label{global-connectedness-theorem-of-Shokurov} Let $(X, B_{X})$
be a log pair with an effective boundary $B_X$. If the divisor
$-(K_{X}+B_{X})$ is nef and big, then the locus $LCS(X, B_{X})$ is
connected.
\end{theorem}

\begin{theorem}
\label{local-connectedness-theorem-of-Shokurov} Let $(X, B_{X})$
be a log pair with an effective boundary. Let $g:X\longrightarrow
Z$ be a contraction. If the  divisor $-(K_{X}+B_{X})$ is  $g$-nef
and $g$-big, then $LCS(X, B_{X})$ is con\-nected in a neighborhood
of each fiber of the contraction $g$.
\end{theorem}

The following result is Theorem 17.4 of \cite{Ko91}.

\begin{theorem}
\label{general-connectedness-theorem-of-Shokurov} Let $g:X\to Z$
be a contraction, where the  varieties $X$ and $Z$ are normal. Let
$D_X=\sum_{i=1}^md_iD_i$ be a $\QQ$-divisor on $X$ such that the
divisor $-(K_X+D_X)$ is $g$-nef and $g$-big. Suppose that
$\codim(g(D_{i})\subset Z)\geq 2$ whenever $d_{i}<0$. Then, for a
log resolution $h:V\to X$ of the log pair $(X, D_X)$, the locus
$\cup_{a_{E}\le -1}E$ is connected in a neighborhood of every
fiber of the morphism $g\circ h$, where $E$ is a divisor on $V$
and the rational number $a_{E}$ is the discrepancy of $E$ with
respect to $(X,D_X)$.
\end{theorem}

\begin{proof}
Let $f=g\circ h$, $A=\sum_{a_{E}>-1}a_E E$, and $B=\sum_{a_{E}\le
-1}-a_E E$. Then
$$\lceil A\rceil-\lfloor B\rfloor=K_{V}-h^{*}(K_{X}+D_{X})+\lbrace-A\rbrace+\lbrace B\rbrace$$
and $R^{1}f_{*}(\OOO_{V}(\lceil A\rceil-\lfloor B\rfloor))=0$ by
Kawamata-Viehweg vanishing. Hence, the map
$$
f_{*}(\OOO_{V}(\lceil A\rceil))\to%
f_{*}(\OOO_{\lfloor B\rfloor}(\lceil A\rceil))
$$
is surjective. Every irreducible component of $\lceil A\rceil$ is
either $h$-ex\-cep\-ti\-onal or the proper transform of some
$D_{j}$ with $d_{j}<0$. Thus $h_{*}(\lceil A\rceil)$ is
$g$-exceptional and $f_{*}(\OOO_{V}(\lceil A\rceil))=\OOO_{Z}$.
Consequently, the map
$$
\OOO_{Z}\to f_{*}(\OOO_{\lfloor B\rfloor}(\lceil A\rceil))
$$
is surjective, which implies the connectedness of $\lfloor
B\rfloor$ in a neighborhood of every fiber of the morphism $f$
because the divisor $\lceil A\rceil$ is effective and has no
common component with $\lfloor B\rfloor$.
\end{proof}

We defined the notions of  centers of canonical singularities and
he set of centers of canonical singularities for movable log
pairs. However, the movability of  boundaries has nothing to do
with all these notions. So we are free to use them for usual log
pairs as well.

The following theorem, frequently referred to as adjunction, leads
us to the bridge between Shokurov's connectedness principle and
N\"other-Fano inequalities.

\begin{theorem}
\label{log-adjunction} Let $(X, B_{X})$ be a log pair with an
effective divisor $B_{X}$, $Z$ an element in $\CCSS(X, B_{X})$,
and $H$  an effective irreducible Cartier divisor on $X$. Suppose
that  both the varieties $X$ and $H$ are smooth at a generic point
of $Z$ and $Z\subset H\not\subset\Supp(B_{X})$. Then, the set
$\LLCCSS(H, B_{X}\vert_{H})$ is not empty.
\end{theorem}

\begin{proof}
Let $f:W\to X$ be a log resolution of $(X, B_{X}+H)$. Put ${\hat
H}=f^{-1}(H)$. Then
$$
K_{W}+{\hat H}= f^{*}(K_{X}+B_{X}+H)+\sum_{E\ne {\hat H}} a(X,
B_{X}+H, E)E
$$
and by our assumption  the subvarieties $Z$ and  $H$ are centers
of the  log canonical singularities of the pair $(X, B_{X}+H)$.
Therefore, applying
Theorem~\xref{general-connectedness-theorem-of-Shokurov} to the
log pullback of $(X, B_{X}+H)$ on $W$, we obtain a divisor  $E\ne
{\hat H}$ on $W$ such that $f(E)=Z$, $a(X, B_{X}, E)\le -1$, and
${\hat H}\cap E\ne\emptyset$. Now the equalities
$$
K_{\hat H}=%
(K_{W}+{\hat H})\vert_{\hat H}= %
f\vert_{\hat H}^{*}(K_{H}+B_{X}\vert_{H})+\sum_{E\ne {\hat H}}
a(X, B_{X}+H, E)E\vert_{\hat H}
$$
imply the claim.
\end{proof}

Our taking N\"other-Fano inequalities into consideration, it is
significant for us to study the singularities of certain movable
log pairs on Fano varieties. It requires us to investigate the
multiplicities of certain movable boundaries or their
self-intersections.

 The following result is Theorem 3.1 of
\cite{Co00}.

\begin{theorem}
\label{theorem-of-Corti} Let $S$ be a smooth surface and
$\MMM_{S}$ an effective movable boundary on the surface $S$.
Suppose that there is a point $o$ in $ \LLCCSS(S,
(1-a_{1})B_{1}+(1-a_{2})B_{2}+\MMM_{S})$, where $a_i$'s are
non-negative rational numbers and $B_{i}$'s are irreducible and
reduced curves on $S$ intersecting normally at the point $o$.
Then, we have
$$
\mult_{o}(\MMM_{S}^{2})\ge\left\{\aligned
&4a_{1}a_{2}\ {\text {if}}\ a_{1}\le 1\ {\text {or}}\ a_{2}\le 1\\
&4(a_{1}+a_{2}-1)\ {\text {if}}\ a_{1}>1\ {\text {and}}\ a_{2}>1.\\
\endaligned
\right.
$$
Furthermore, the inequality is strict if the singularities of the
log pair $(S, (1-a_{1})B_{1}+(1-a_{2})B_{2}+\MMM_{S})$ are not log
canonical in a neighborhood of the point $o$.
\end{theorem}

\begin{proof}
Let $D=(1-a_{1})B_{1}+(1-a_{2})B_{2}+\MMM_{S}$ and $f:S'\to S$ be
a birational morphism such that the surface $S'$ is smooth. We
consider
$$
K_{S'}+f^{-1}(D)=
E_{i})E_{i},
$$
where $E_{i}$ is an $f$-exceptional curve. We suppose that $a(S,
D, E_{1})\le -1$ and the curve $E_1$ is contracted to the point
$o$. Then the birational morphism $f$ is a composition of $k$ blow
ups at smooth points.

\textbf{\itshape Claim 1.}  The statement is true when $a_{1}>1$
and $a_{2}>1$ if the statement holds when $a_{1}\le 1$ or
$a_{2}\le 1$.

Define the numbers $a(S, E_{i})$, $m(S, \MMM_{S}, E_{i})$ and
$m(S, B_{j}, E_{i})$ as follows;
$$
\aligned
&\sum_{i=1}^{k}a(S, E_{i})E_{i}= K_{S'}-f^{*}(K_{S}),\\
&\sum_{i=1}^{k}m(S, \MMM_{S}, E_{i})E_{i}=
f^{-1}(\MMM_{S})-f^{*}(\MMM_{S}),\\
&\sum_{i=1}^{k}m(S, B_{j}, E_{i})E_{i}=
f^{-1}(B_{j})-f^{*}(B_{j}).\\
\endaligned
$$
We then observe that the equality
$$
a(S, D, E_{i})=a(S, E_{i})-m(S, \MMM_{S}, E_{i})+%
m(S, B_{1}, E_{i})(a_{1}-1)+m(S, B_{2}, E_{i})(a_{2}-1)
$$
holds. We may assume that %
$m(S, B_{1}, E_{1})\geq m(S, B_{2}, E_{1})$. %
Then,
$$
-1\geq a(S, D, E_{1})\geq %
a(S, E_{1})-m(S, \MMM_{S}, E_{1})+m(S, B_{2},
E_{1})(a_{1}+a_{2}-2)
$$
and hence $o\in \LLCCSS(S, (2-a_{1}-a_{2})B_{2}+\MMM_{S})$. %
Because the log pair  %
$(S, (2-a_{1}-a_{2})B_{2}+\MMM_{S})$ satisfies our assumption, we
obtain $\mult_{o}(\MMM_{S}^{2})\geq 4(a_{1}+a_{2}-1)$.

\textbf{\itshape Claim 2.}  The statement holds when $a_{1}\le 1$
or $a_{2}\le 1$.

We may assume that $a_{1}\le 1$. Let $h:T\to S$ be the blow up at
the point $o$ and $E$ be an exceptional curve of $h$. Then $f$
factors through $h$ such that $f=g\circ h$ for some birational
morphism $g:S'\to T$ which is a composition of $k-1$ blow ups at
smooth points. Then
$$
K_{T}+(1-a_{1}){\bar B}_{1}+%
(1-a_{2}){\bar B}_{2}+(1-a_{1}-a_{2}+m)E+\MMM_{T}= %
h^{*}(K_{S}+D),
$$
where ${\bar B}_{j}=h^{-1}(B_{j})$, $m=\mult_{o}(\MMM_{S})$, and
$\MMM_{T}=h^{-1}(\MMM_{S})$.
\par
We are to use the induction on $k$. In the case $k=1$, we have
$S'=T$, $E_{1}=E$, and $a(S, D, E_{1})=a_{1}+a_{2}-m-1\le -1$.
Thus
$$
\mult_{o}(\MMM_{S}^{2})\geq m^{2}\geq (a_{1}+a_{2})^{2}\geq
4a_{1}a_{2}
$$
and we are done.

We therefore suppose that $k>1$ and $g(E_{1})$ is a point $p\in
E$. We see
$$
p\in \LLCCSS(T, (1-a_{1}){\bar B}_{1}+%
(1-a_{2}){\bar B}_{2} +(1-a_{1}-a_{2}+m)E+\MMM_{T}).
$$
There are three possible cases:
$p\in E\cap {\bar B}_{1}$, %
$p\in E\cap {\bar B}_{2}$, and %
$p\not\in {\bar B}_{1}\cup {\bar B}_{2}$. %
By  the induction hypothesis,   the statement holds for the log
pair
$$
(T, (1-a_{1}){\bar B}_{1}+(1-a_{1}-a_{2}+m)E+\MMM_{T})
$$
in the case $p\in E\cap {\bar B}_{1}$, for the log pair
$$
(T, (1-a_{2}){\bar B}_{2}+(1-a_{1}-a_{2}+m)E+\MMM_{T})
$$
in the case $p\in E\cap {\bar B}_{2}$, and  for the log pair %
$$(T, (1-a_{1}-a_{2}+m)E+\MMM_{T})$$ in the case $p\not\in {\bar
B}_{1}\cup {\bar B}_{2}$ because all conditions of the theorem are
satisfied in these cases and the morphism $g$ consists of $k-1$
blow ups at smooth points. Also we have
$$\mult_{o}(\MMM_{S}^{2})\geq m^{2}+\mult_{p}(\MMM_{T}^{2}).$$
\par
In the case $p\in E\cap {\bar B}_{1}$, we obtain
$$
\mult_{o}(\MMM_{S}^{2})\geq
m^{2}+4a_{1}(a_{1}+a_{2}-m)=(2a_{1}-m)^{2}+4a_{1}a_{2}\geq
4a_{1}a_{2}.
$$
\par
Consider the case $p\in E\cap {\bar B}_{2}$. If either $a_{2}\le
1$ or $a_{1}+a_{2}-m\le 1$, then we can proceed as in the previous
case. If not, then we have
$$
\mult_{o}(\MMM_{S}^{2})\geq m^{2}+4(a_{1}+2a_{2}-m-1)>4a_{2}\geq
4a_{1}a_{2}.
$$
\par
If $p\not\in {\bar B}_{1}\cup {\bar B}_{2}$, then we obtain
$$
\mult_{o}(\MMM_{S}^{2})\geq
m^{2}+4(a_{1}+a_{2}-m)>m^{2}+4a_{1}(a_{1}+a_{2}-m)\geq
4a_{1}a_{2},
$$
which completes the proof.
\end{proof}

Instead of Theorem~\ref{theorem-of-Corti}, the following
simplified version, which is a special case of Theorem 2.1 in
\cite{dFEM03}, is  more often applied.

\begin{theorem}
\label{theorem-of-Corti-simple-form} Let $S$ be a smooth surface,
$o$ a point on $S$, and $\MMM_{S}$ an effective movable boundary
on $S$ such that $o\in \LLCCSS(S, \MMM_{S})$. Then
$\mult_{o}(\MMM_{S}^{2})\geq 4$.  Moreover, if the equality holds,
then $\mult_o(\MMM_{S})=2$.
\end{theorem}

Even though Theorems~\ref{theorem-of-Corti} and
~\ref{theorem-of-Corti-simple-form} are results on surfaces, they
can be applied to 3-folds via Theorem~\ref{log-adjunction}. The
following result is Corollary 7.3 of \cite{Pu00}, which holds even
over fields of positive characteristic and implicitly goes back to
the classical paper \cite{IsMa71}.

\begin{theorem}
\label{theorem-of-Iskovskikh} Let $X$ be a smooth 3-fold and
$\MMM_{X}$  an effective movable boundary on $X$. Suppose that a
point $o$ belongs to $\CCSS(X, \MMM_{X})$. Then the inequality
$\mult_{o}(\MMM_{X}^{2})\geq 4$ holds and the equality holds only
when $\mult_{o}(\MMM_{X})=2$.
\end{theorem}

\begin{proof} Let $H$ be a general very ample divisor on $X$
containing $o$. Then the point $o$ is a center of log canonical
singularities of the log pair $(H, \MMM_{X}\vert_{H})$ by
Theorem~\xref{log-adjunction}. On the other hand,
$$
\mult_{o}(\MMM_{X}^{2})=\mult_{o}((\MMM_{X}\vert_{H})^{2})
$$
and $\mult_{o}(\MMM_{X})=\mult_{o}(\MMM_{X}\vert_{H})$. Hence, the
claim follows from Theorem~\xref{theorem-of-Corti-simple-form}.
\end{proof}

As a matter of fact, Theorem~\xref{theorem-of-Iskovskikh} can be
proved in a more geometric way.

\begin{lemma}
\label{lemma-on-terminal-extraction} Let $X$ be a smooth 3-fold
and $\MMM_{X}$ an effective movable boundary on $X$. Suppose that
the log pair $(X, \MMM_{X})$ has canonical singularities and
$\CCSS(X, \MMM_{X})$ contains a point $o$. Then there is a
birational morphism $f:V\to X$ such that $V$ has $\QQ$-factorial
terminal singularities, $f$ contracts exactly one exceptional
divisor $E$ to the point $o$, and
$$
K_{V}+\MMM_{V}=f^{*}(K_{X}+\MMM_{X}),
$$
where $\MMM_{V}=f^{-1}(\MMM_{X})$.
\end{lemma}

\begin{proof}
Because the log pair $(X, \MMM_X)$ has at worst canonical
singularities, there are finitely many divisorial discrete
valuations $\nu$ of the field of rational functions of $X$ whose
centers on $X$ are the point $o$ and whose discrepancies
$a(X,\MMM_X,\nu)$ are non-positive.
 Therefore, we may consider a birational morphism
$g:W\to X$ such that the 3-fold $W$ is smooth, $g$ contracts $k$
divisors,
$$
K_{W}+\MMM_{W}= g^{*}(K_{X}+\MMM_{X})+\sum_{i=1}^{k}a_{i}E_{i},
$$
the movable log pair $(W, \MMM_{W})$ has canonical singularities,
and the set $\CCSS(W, \MMM_{W})$ does not contain subvarieties of
$\cup_{i=1}^{k}E_{i}$, where $\MMM_{W}=g^{-1}(\MMM_{X})$,
$g(E_{i})=o$, and $a_{i}\in \QQ$. Applying the relative version of
Log Minimal Model Program (see \cite{KMM}) to the movable log pair
$(W, \MMM_{W})$ over $X$, we may assume that $W$ has
$\QQ$-factorial terminal
 singularities and
$$
K_{W}+\MMM_{W}= g^{*}(K_{X}+\MMM_{X})
$$
because of the canonicity of $(X, \MMM_{X})$. Applying the
relative Minimal Model Program to $W$ over the variety $X$, we get
the necessary 3-fold and the birational morphism.
\end{proof}

The following result was conjectured in \cite{Co95} and proved in
\cite{Kaw01}.

\begin{theorem}
\label{theorem-of-Kawakita} Let $X$ be a smooth 3-fold and $f:V\to
X$ be a birational morphism of a 3-fold $V$ with $\QQ$-factorial
terminal
 singularities. Suppose that the morphism $f$
contracts exactly one exceptional divisor $E$ and contracts it to
a point $o$. Then the morphism $f$ is the weighted blow up at the
point $o$ with weights $(1,n_1,n_2)$ in suitable local coordinates
on $X$, where the natural numbers $n_1$ and $n_2$ are coprime.
\end{theorem}

With Theorem~\xref{theorem-of-Kawakita},
Theorem~\xref{theorem-of-Iskovskikh} was proved in \cite{Co95} in
the following way, which explains the geometrical nature of the
inequality in Theorem~\xref{theorem-of-Iskovskikh}.

\begin{proposition}
\label{proposition-on-weighted-intersections} Let $X$ be a smooth
3-fold with an effective movable boundary $\MMM_X$ on $X$. Suppose
that $\CCSS(X, \MMM_{X})$ contains a point $o$. Let $f:V\to X$ be
the weighted blow up at the point $o$ with weights $(1,n_1,n_2)$
in suitable local coordinates on $X$ such that
$$
K_{V}+\MMM_{V}=f^{*}(K_{X}+\MMM_{X}),
$$
where natural numbers $n_1$ and $n_2$ are coprime and
$\MMM_{V}=f^{-1}(\MMM_{X})$. Then
$$
\mult_{o}(\MMM_{X}^{2})\geq {\frac {(n_1+n_2)^{2}}
{n_1n_2}}=4+{\frac {(n_1-n_2)^{2}} {n_1n_2}}\geq 4.
$$
Moreover, if $n_1=n_2$, then $f$ is the regular blow up at the
point $o$ and $\mult_{o}(\MMM_{X})=2$.
\end{proposition}

\begin{proof} Let $E\subset V$ be the $f$-ex\-cep\-ti\-o\-nal
divisor. Then
$$
K_{V}= f^{*}(K_{X})+(n_1+n_2)E
$$
and $\MMM_{V}=f^{*}(\MMM_{X})-mE$ for some $m\in \QQ_{>0}$. Thus,
$m=n_1+n_2$ and
$$
\mult_{o}(\MMM_{X}^{2})\geq m^{2}E^{3}={\frac {(n_1+n_2)^{2}}
{n_1n_2}}.
$$
\end{proof}

The following application of Theorem~\xref{log-adjunction} is
Theorem 3.10 in \cite{Co00}.

\begin{theorem}
\label{theorem-about-double-point} Let $X$ be a 3-fold with a
simple double point $o$ and $B_{X}$ an effective boundary on $X$
such that $o\in \CCSS(X, B_{X})$. Then the inequality
$\mult_{o}(B_{X})\geq 1$ holds.
\end{theorem}

\begin{proof}
Let $f:W\to X$ be the blow up at the point $o$ and $E$ be the
$f$-exceptional divisor. Then
$$
K_{W}+B_{W}= f^{*}(K_{X}+B_{X})+(1-\mult_{o}(B_{X}))E,
$$
where $B_{W}=f^{-1}(B_{X})$. Suppose that $\mult_{o}(B_{X})<1$.
Then, there is a center $Z\in \CCSS(W, B_{W})$ that is contained
in $E$, and hence
$$
\LLCCSS(E, B_{W}\vert_{E})\ne\emptyset
$$
by Theorem~\xref{log-adjunction}. But it is impossible because of
Lemma~\xref{lemma-on-quadric}.
\end{proof}


\section{Birational
super-rigidity.}\label{section-birational-super-rigidity}

The goal of this section is to prove Theorem~A.

Let $\pi:X\to \PP^3$ be a $\QQ$-factorial double cover ramified
along a nodal sextic $S\subset \PP^3$. We then see that
$\Pic(X)\cong \ZZ K_{X}$, $ -K_{X}\sim \pi^{*}(\OOO_{\PP^3}(1)), $
and $-K_{X}^{3}=2$. Consider an arbitrary movable boundary
$\MMM_{X}$ on the $3$-fold $X$ such that the divisor
$-(K_{X}+\MMM_{X})$ is ample. To prove Theorem~A we must show that
$\CCSS(X, \MMM_{X})=\emptyset$ and then apply
Theorem~\xref{theorem-Nother-Fano-inequality}.

We suppose that $\CCSS(X, \MMM_{X})\ne\emptyset$. In what follows,
we will derive a contradiction.

\begin{lemma}
\label{lemma-on-smooth-points} Smooth points of the  $3$-fold $X$
are not contained in $\CCSS(X, \MMM_{X})$.
\end{lemma}

\begin{proof}
Suppose that $\CCSS(X, \MMM_{X})$ has a smooth point $o$ on $X$.
Let $H$ be a general enough divisor in the linear system
$|-K_{X}|$ passing through the point $o$. We then obtain

$$
2=H\cdot K_{X}^{2}>H\cdot \MMM_{X}^{2}\geq
\mult_{o}(\MMM_{X}^{2})\geq 4
$$
from Theorem~\xref{theorem-of-Iskovskikh}, which is absurd.
\end{proof}

\begin{lemma}
\label{lemma-on-singular-points} Singular points of the $3$-fold
$X$ are not contained in $\CCSS(X, \MMM_{X})$.
\end{lemma}

\begin{proof}
If $\CCSS(X, \MMM_{X})$ contains a singular point $o$ on $X$, then
Theorem~\xref{theorem-about-double-point} gives us
$$
2=H\cdot K_{X}^{2}>H\cdot \MMM_{X}^{2}\geq
2mult^{2}_{o}(\MMM_{X})\geq 2,
$$
where $H$ is a general enough divisor in $|-K_{X}|$ passing
through the point $o$. It is absurd.
\end{proof}
Lemmas~\xref{lemma-on-smooth-points}
and~\xref{lemma-on-singular-points} together show that any element
of the set  $\CCSS(X, \MMM_{X})$ cannot be a point of $X$.
Therefore, it must contain a curve $C\subset X$. To complete the
proof of Theorem~A it is enough to show that the set $\CCSS(X,
\MMM_{X})$ cannot contain a curve.

\begin{lemma}
\label{lemma-on-degree-of-curve} The intersection number
$-K_{X}\cdot C$ is $1$.
\end{lemma}

\begin{proof}
Let $H$ be a general enough divisor in the anticanonical linear
system $|-K_{X}|$. Then
$$
2=H\cdot K_{X}^{2}>H\cdot \MMM_{X}^{2}\geq
\mult_{C}(\MMM_{X}^{2})H\cdot C\geq -K_{X}\cdot C,
$$
which implies $-K_{X}\cdot C=1$.
\end{proof}

\begin{corollary}
\label{corollary-on-image-of-curve} The curve  $\pi(C)\subset
\PP^3$ is a line and $C\cong \PP^1$.
\end{corollary}

\begin{lemma}
\label{lemma-on-curve-in-smooth-locus} The curve $C$ is not
contained in the smooth locus of the $3$-fold $X$.
\end{lemma}

\begin{proof}
Suppose that the curve $C$ lies on the smooth locus of the
$3$-fold $X$. Let $f:W\to X$ be the blow up along the curve $C$
and $E$ be the $f$-ex\-cep\-tional divisor. We then get
$\mult_{C}(\MMM_{X})\geq 1$ and
$$
\MMM_{W}=f^{-1}(\MMM_{X})= f^*(\MMM_{X})-\mult_{C}(\MMM_{X})E.
$$
The linear system $|-K_{W}|=|f^{*}(-K_{X})-E|$ has just one base
curve ${\tilde C}$ such that
$$
\pi\circ f(\tilde C)=\pi(C)\subset\PP^3.
$$
We see that ${\tilde C}\subset E$ if and only if $\pi(C)\subset
S$.

Let $H=f^{*}(-K_{X})$. Then the divisor $3H-E$ has non-negative
intersection with all the curves on $W$ possibly except $\tilde
C$. We are to show that the divisor $3H-E$ is nef. We obtain
$(3H-E)\cdot {\tilde C}=0$ unless ${\tilde C}$ is contained in
$E$. Therefore, we suppose that the curve ${\tilde C}$ is
contained in $E$.

The normal bundle $\NNN_{X/C}$ of the curve $C\cong \PP^1$ on the
$3$-fold $X$ splits into
$$
\NNN_{X/C}\cong\OOO_{C}(a)\oplus \OOO_{C}(b)
$$
for some integers $a\geq b$. The exact sequence
$$
0\to \TTT_{C}\to \TTT_{X}\vert_{C}\to \NNN_{X/C}\to 0
$$
shows $\deg(\NNN_{X/C})=a+b=-K_{X}\cdot C+2g(C)-2=-1$.

On the other hand, the curve $C$ is contained in the smooth locus
of the proper transform ${\hat S}\cong S$ of the sextic $S\subset
\PP^3$. The exact sequence
$$
0\to \NNN_{{\hat S}/C}\to \NNN_{X/C}\to \NNN_{X/{\hat S}}\vert_{C}\to 0%
$$
and $\NNN_{{\hat S}/C}\cong \OOO_{C}(-4)$ imply $b\geq -4$. In
particular, $a-b\le 7$.

Let $s_{\infty}$ be the exceptional section of the ruled surface
$f\vert_{E}:E\to C$. Because $E^{3}=-\deg(\NNN_{X/C})=1$ and
$-K_X\cdot C=1$, we obtain
$$
(3H-E)\cdot s_{\infty}={\frac {7+b-a} {2}}\geq 0,
$$
which implies that the divisor $3H-E$ is nef.

Because $3H-E$ is nef, we get $(3H-E)\cdot \MMM_{W}^{2}\geq 0$,
but
$$
(3H-E)\cdot \MMM_{W}^{2}=6r^{2}-4mult^{2}_{C}(\MMM_{X})-2r
\mult_{C}(\MMM_{X})<0,
$$
where $r\in \QQ\cap (0,1)$ such that $\MMM_{X}\sim_\QQ -r K_{X}$.
\end{proof}

\begin{corollary}
\label{corollary-on-singular-points-on-curve} The curve $C$
contains a simple double point of the $3$-fold $X$.
\end{corollary}

\begin{lemma}
\label{lemma-on-good-curve} The line $\pi(C)$ is contained in the
sextic surface $S$.
\end{lemma}

\begin{proof}
Suppose $\pi(C)\not\subset S$. Let $\HHH$ be the  linear subsystem
in $|-K_{X}|$ of surfaces containing the curve $C$. The base locus
of $\HHH$ consists of the curve $C$ and the curve ${\widetilde C}$
such that $\pi(C)=\pi({\widetilde C})$.  Choose a general enough
surface $D$ in the pencil $\HHH$.  The restriction
$\MMM_{X}\vert_{D}$ is not movable, but
$$
\MMM_{X}\vert_{D}=\mult_{C}(\MMM_{X})C+\mult_{\widetilde
C}(\MMM_{X}){\tilde C}+\mathcal{R}_{D},
$$
where $\mathcal{R}_{D}$ is a movable boundary. The surface $D$ is
smooth outside of the singular points $p_{i}$ of the $3$-fold $X$
which are contained in the curve $C$. Moreover, each point $p_{i}$
is a simple double point on the surface $D$. Thus, on the surface
$D$, we have
$$C^{2}={\tilde C}^{2}=-2+{\frac {k} {2}},$$
where $k$ is the number of the points $p_{i}$ on $C$. Hence, we
obtain $C^2={\widetilde C}^{2}<0$ on the surface $D$ because
$k\leq 3$. Immediately, the inequality
$$
(1-\mult_{\tilde C}(\MMM_{X})){\widetilde C}^{2}\geq
(\mult_{C}(\MMM_{X})-1)C\cdot {\widetilde C}+\mathcal{R}_{D}\cdot
{\widetilde C}\geq 0
$$
follows, which implies $\mult_{\tilde C}(\MMM_{X})\geq 1$.
Therefore, for a general member $H\in |-K_{X}|$ we have a
contradiction
$$
2=H\cdot K_{X}^{2}>H\cdot \MMM_{X}^{2}\geq
mult^2_{C}(\MMM_{X})H\cdot C+mult^2_{\tilde C}(\MMM_{X})H\cdot
{\tilde C}\geq 2.
$$
\end{proof}

\begin{lemma}
\label{lemma-on-bad-curve} The line $\pi(C)$ is not contained in
the sextic surface $S$.
\end{lemma}

\begin{proof}
Suppose $\pi(C)\subset S$. Let $p$ be a general enough point on
the curve $C$ and $L\subset \PP^3$ be a general line tangent to
$S$ at the point $\pi(p)$. Then the proper transform ${\tilde
L}\subset X$ of $L$ is an irreducible curve which is singular at
the point $p$. By construction, the curve $L$ is not contained in
the base locus of the components of the movable boundary
$\MMM_{X}$. Thus, we obtain contradictory inequalities
$$
2>{\tilde L}\cdot\MMM_{X}\ge \mult_{p}{\tilde
L}\mult_{p}(\MMM_{X})\ge 2\mult_{C}(\MMM_{X})\ge 2.
$$
\end{proof}

We have shown that the set $\CCSS(X, \MMM_{X})$ is empty. Now, we
can immediately  obtain Theorem~A from
Theorem~\xref{theorem-Nother-Fano-inequality}.


\section{$\QQ$-factoriality.}

In this section we  study the $\QQ$-factoriality on double covers
of $\PP^3$ ramified along  nodal sextics   and prove Theorem~B.

The $\QQ$-factoriality depends both on local types of
singularities and on their global position. In the case of Fano
$3$-folds, the $\QQ$-factoriality is equivalent to the global
topological condition
$$
\rank(H^2(X,\ZZ))=\rank(H_4(X,\ZZ)).
$$
 In the case of the double solids, the condition means
the 4th integral homology group of $X$ generated by the class of
the full back of a hyperplane in $\PP^3$ via the covering
morphism.

Using the method in \cite{Cl83}, we study the $\QQ$-factoriality
on a double cover $X$ of $\PP^3$ ramified along a sextic $S$. As
before, we assume that $X$ has only simple double points. Note
that $\Pic(X)\cong H^2(X,\ZZ)$ when $X$ has at worst rational
singularities.

For us in order to see whether a double solid $X$ is
$\QQ$-factorial, the main job is to compute the rank of the group
$H_4(X, \ZZ)$. Indeed, the double solid $X$ is $\QQ$-factorial if
and only if $\rank(H_4(X, \ZZ))=1$ because $\rank H^2(X,\ZZ)=1$.
The paper \cite{Cl83} gives us a method to compute it by studying
the number of singularities of $S$, their position in $\PP^3$, and
the sheaf $\III\otimes\OOO_{\PP^3}(5)$, where $\III$ is the ideal
sheaf of the set $\Sigma$ of singular points of $S$ in
$\mathbb{P}^3$. The following result was proved in \cite{Cl83}
(see also \cite{Di90} and \cite{Cy01}).

\begin{theorem}
\label{theorem-of-Clemens} Under the same notation, we have
$$
\rank (H_4(X,\ZZ))= \# (\Sigma)-I +1,
$$
where $I$ is the number of independent conditions which vanishing
on $\Sigma$ imposes on homogeneous forms of degree $5$ on $\PP^3$.
\end{theorem}

We define the \emph{defect} of $X$ to be the non-negative integer
$\# (\Sigma)-I$. Then we can restate the $\QQ$-factoriality as
follows:

\begin{corollary}
The double cover $X$ is $\QQ$-factorial if and only if the defect
of $X$ is $0$.
\end{corollary}

On the other hand, from the exact sequence
$$
0\to \III\otimes\OOO_{\PP^3}(5)\to \OOO_{\PP^3}(5)\to
\bigoplus_{p\in\Sigma}\CC\to 0
$$
we obtain a long exact sequence
$$
0\to H^0(\PP^3, \III\otimes\OOO_{\PP^3}(5))\to
H^0(\PP^3,\OOO_{\PP^3}(5))\to
H^0(\PP^3,\bigoplus_{p\in\Sigma}\CC)\to H^1(\PP^3,
\III\otimes\OOO_{\PP^3}(5))\to 0,
$$
which tells us
\[\mbox{defect of } X=\dim(H^1(\PP^3,
\mathcal{I}\otimes\mathcal{O}_{\PP^3}(5))).\]

An immediate application of this method is the second part of
Theorem~B. Since $\dim(H^0(\PP^3,\mathcal{O}_{\PP^3}(5)))=56$, the
defect of $X$ is positive if $\#(\Sigma)\geq 57$.

We can easily observe that if $\#(\Sigma)\leq 6$, then the
sequence
$$
0\to H^0(\PP^3, \mathcal{I}\otimes\mathcal{O}_{\PP^3}(5))\to
H^0(\PP^3,\mathcal{O}_{\PP^3}(5))\to
H^0(\PP^3,\bigoplus_{p\in\Sigma}\CC)\to 0,
$$
is exact regardless of their position. Therefore, when
$\#(\Sigma)\leq 6$ the defect of $X$ is trivially $0$,
\emph{i.e.}, the sextic double solid $X$ is $\QQ$-factorial. As a
matter of fact, we can go farther. As Theorem~B states, if the
surface $S$ has at most $14$ nodes, then the $3$-fold $X$ is
$\QQ$-factorial regardless of their position. In what follows, we
prove the first part of Theorem~B.

\begin{definition}
We say that a set of points $\Gamma$ on $\PP^3$ is on sextic-node
position if no $5k+1$ points of $\Gamma$ can lie on a curve of
degree $k$ in $\PP^3$ for every positive integer $k$.
\end{definition}
\begin{lemma}
\label{lemma-on-number-of-double-points-on-curves} Let $\Sigma$ be
the set of singular points of the sextic $S$. Then the set
$\Sigma$ is on sextic-node position.
\end{lemma}

\begin{proof}
Suppose that the surface $S$ is defined by a homogeneous
polynomial equation $F(x_0,x_1,x_2,x_3)=0$ of degree six. We
consider the linear system
\[\mathcal{L}:=\left|\sum_{i=0}^{3}\lambda_i\frac{\partial F}{\partial x_i}=0\right|.\]
The base locus of the linear system $\mathcal{L}$ is exactly the
singular locus of the surface $S$. A curve of degree $k$ in
$\PP^3$ intersects a generic member of the linear system
$\mathcal{L}$ exactly $5k$ times since
$\mathcal{L}\subset|\mathcal{O}_{\PP^3}(5)|$. Therefore, the set
$\Sigma$ is on sextic-node position.
\end{proof}

For convenience, we state an elementary lemma.
\begin{lemma}
\label{lemma-of-excluding-three-points} Let $\Gamma=\{q_1,\cdots,
q_s\}$ be a set of $s\geq 4$ points in $\PP^3$. For a given point
$q\not\in\Gamma$, there is a hyperplane $H$ which contains at
least three points of $\Gamma$ but not the point $q$ unless all
the points $q, q_1, \cdots, q_s$ lie on a single hyperplane.
\end{lemma}
\begin{proof}
Because not all the points $q, q_1,\cdots, q_s$ lie on a single
hyperplane, we may assume there are two distinct hyperplane $H_1$
and $H_2$ such that $H_1\cup H_2$ contains the point $q$ and four
points, say $q_1, q_2, q_3$, and $q_4$, of $\Gamma$, $q_1\in
H_1\setminus H_2$, and $q_2\in H_2\setminus H_1$. Then one of the
hyperplanes generated by $\{q_1, q_2, q_3\}$ and $\{q_1, q_2,
q_4\}$ must not pass through the point $q$; otherwise all of the
five points $q, q_1,\cdots, q_4$ would be on a single  hyperplane.
\end{proof}

Also, the following result of \cite{Bes83} is useful.

\begin{theorem}
\label{theorem-of-Bese} Let $\pi:Y\to\PP^2$ be the blow up at
points $p_1, \cdots, p_s$ on $\PP^2$. Then the linear system
$|\pi^{*}(\OOO_{\PP^2}(d))-\sum_{i=1}^{s}E_{i}|$ is
base-point-free for all $s\leq \frac{1}{3}(h^0(\PP^2,
\OOO_{\PP^2}(d+3) )-5)$, where $d\geq 3$ and $E_i=\pi^{-1}(p_i)$,
if at most $k(d+3-k)-2$ of the points $p_i$ lie on a curve of
degree $1\leq k\leq \frac{1}{2}(d+3)$.
\end{theorem}

Theorem~\ref{theorem-of-Clemens} tells us that the first part of
Theorem~B immediately follows from the lemma below.

\begin{lemma}
\label{lemma-on-blow-ups-of-points} Let $\gamma:V\to\PP^3$ be the
blow up at $k$ different points $\Gamma=\{p_{1},\ldots,p_{k}\}$
and $p$ be a point in $V\setminus\cup_{i=1}^kE_i$ such that the
set $\Gamma\cup\{\gamma(p)\}$ is on sextic-node position, where
$E_{i}=\gamma^{-1}(p_{i})$. If $k\leq 13$, then the linear system
$|\gamma^{*}(\OOO_{\PP^3}(5))-\sum_{i=1}^{k}E_{i}|$ is
base-point-free at the point $p$.
\end{lemma}
\begin{proof}
It is enough to find a quintic hypersurface in $\PP^3$ that passes
through all the points of $\Gamma$ but not the point
$q:=\gamma(p)$.  We may assume that $k=13$. Let $r$ be the maximal
number of points of $\Gamma$ that belong to a single hyperplane of
$\PP^3$ together with the point $q$. Note that $2\leq r\leq 13$.
Without loss of generality, we may also assume that the first $r$
points of $\Gamma$, \emph{i.e.}, $p_1,\cdots, p_r$, are contained
in a hyperplane $H$  together with the point $q$.

We prove the statement case by case.

\emph{Case 1.} (r=2)

We divide the set $\Gamma$ into five subsets of $\Gamma$ such that
each subset contains exactly three points of $\Gamma$ and the
union of all the five subsets is $\Gamma$. Because $r=2$, the
hyperplane generated by each subset cannot contains the point $q$.
The product of these five hyperplanes is what we want.

Before we proceed, we  note that the points $q$ and $p_1, \cdots,
p_r$ do not lie on a single line. If they do, then the hyperplane
$H$ must contain more than $r$ points of $\Gamma$. \vspace{5mm}

\emph{Case 2.} (r=3)

By Lemma~\ref{lemma-of-excluding-three-points}, we can find three
points of $\Gamma$ outside of $H$ such that generate a hyperplane
not passing though the point $q$. Since $r=3$, we can repeat this
procedure two more times with the remaining points of $\Gamma$ in
the outside of $H$. Only one point, say $p_{13}$, then remains in
the  outside of $H$. Because the four points, $q$, $p_1, p_2,p_3$,
cannot lie on a line, there is a quadric hypersurface passing
through the points $p_1,p_2, p_3, p_{13}$, but not the point $q$.

\emph{Case 3.} (r=4)

As in the previous case, we can find two hyperplanes  which
together contains six points of $\Gamma$ in the outside of $H$ but
not $q$. We then have three remaining points of $\Gamma$ in the
outside of $H$. There is a line passing though two points, say
$p_1, p_2$, of $p_1,\cdots, p_4$, but not the point $q$. Then
these two points together with one of the remaining points in the
outside of $H$ generate a hyperplane not containing the point $q$.
Now, we have four points, two of them are on $H$ and the others
not on $H$. Obviously, these four points belong to a quadric
hypersurface not passing through the point $q$. Therefore, the
product of the quadric hypersurface and the hyperplanes gives us a
quintic hypersurface that we are looking for.

\emph{Case 4.} (r=5)

First of all, by Lemma~\ref{lemma-of-excluding-three-points}, we
find a hyperplane which contains three points, say $p_6, p_7,
p_8$, of $\Gamma$ in the outside of $H$ but not the point $q$. We
split the case into two subcases.

\emph{Subcase 4.1.} When four points of $\Gamma$ on $H$ together
with the point $q$ lie on a line.

Assume that the points $q$ and $p_1,\cdots, p_4$ lie on a single
line. The hyperplane generated by the points $p_4$, $p_5$, and
$p_9$ cannot contains $q$. The hyperplane generated by $\{p_3,
p_{10}, p_{11}\}$ cannot pass through the point $q$; otherwise the
number $r$ would be bigger than five. By the same reason, we can
find a hyperplane which contains $\{p_2, p_{12}, p_{13}\}$ but not
the point  $q$. Choose a hyperplane which passes through the point
$p_1$ but not the point $q$. Then we are done.

\emph{Subcase 4.2.} When no four points of $\Gamma$ on $H$ lie on
a line together with the point $q$.

In this case, two pairs of  points of $\Gamma$ on $H$ give two
lines which do not contain the point $q$. Therefore, we can find a
quadric hypersurface containing six points of $\Gamma$, four from
$H$ and two from $\Gamma\setminus (H\cup\{p_6,p_7, p_8\})$, but
not the point $q$. Furthermore, because the number $r$ is five  we
may choose the two points from $\Gamma\setminus (H\cup\{p_6,p_7,
p_8\})$ so that the other three points in the outside of $H$
cannot belong to a single line together with the point $q$. We
then have four points which we have not covered yet, three points,
say $p_{11}, p_{12}, p_{13}$, from the outside of $H$, and one
point, say $p_1$, on $H$. Because the points $p_{11}, p_{12},
p_{13}$ and $q$ do not lie on a line, we can easily find a quadric
hypersurface passing through all the four points but not the point
$q$.

\emph{Case 5.} (r=6)

Again, by Lemma~\ref{lemma-of-excluding-three-points}, we find a
hyperplane which contains three points, say $p_7, p_8, p_9$, of
$\Gamma$ in the outside of $H$ but not the point $q$. By the
sextic-node position condition, we can find two lines on $H$ which
together contain four points of $\Gamma$ on $H$ but not the point
$q$. They give us a quadric hypersurface in $\PP^3$ which pass
though six points of $\Gamma\setminus\{p_7, p_8, p_9\}$. Among
these six points, two points are from the outside of $H$ and the
others from $H$. Therefore, we have four points that have not been
yet covered. Because two of them are in the outside of $H$, we can
easily find a quadric hypersurface which contains these four
points but not the point $q$.

\emph{Case 6.} (r=7)

In this case, we can find three pairs of points of $\Gamma$ on $H$
such that each pair gives us a line not passing through the point
$q$. It implies that we can construct a cubic hypersurface which
passes through six points of $\Gamma$ on $H$ and three points of
$\Gamma$ in the outside of $H$ but not the point  $q$. Moreover,
we may assume that the remaining three points in the outside of
$H$ do not lie on a single line together with the point $q$ due to
the sextic-node position condition. It is easy to find a quadric
hypersurface containing the remaining points of $\Gamma$ but not
$q$. So we are done.

\emph{Case 7.} (r=8 or 9)

We can find four pairs of points of $\Gamma$ on $H$ such that each
pair gives us a line not passing through the point $q$. From this
fact, we easily obtain a quartic hypersurface passing eight points
of $\Gamma$ on $H$ and four points of $\Gamma$ outside of $H$ but
not the point $q$. We then have only one point of $\Gamma$ that is
not covered. Just take a hyperplane passing through this point but
not the point $q$, and we are done.

\emph{Case 8.} (r=10)

Because of the sextic-node position condition, we can find three
pairs, say $\{p_1, p_2\}$, $\{p_3, p_4\}$, $\{p_5, p_6\}$ of
points of $\Gamma$ on $H$ such that each pair gives us a line not
passing through the point $q$ and, in addition, no three of the
points $p_7, p_8, p_9, p_{10}$ cannot lie on a line passing
through point $q$. This shows there is a quintic hypersurface
which passes through $\Gamma$ but not the point $q$.

\emph{Case 9.} (r=11)

We have eleven points of $\Gamma$ on $H$ and two points, $p_{12},
p_{13}$ of $\Gamma$ in the outside of $H$. We can find a quintic
curve $C$ on $H$ which passes through the eleven points on $H$ but
not the point $q$ by Theorem~\ref{theorem-of-Bese}. Note that the
support of the curve $C$ is not a line because of the sextic-node
position condition. A generic hyperplane passing through $p_{12},
p_{13}$ meets $C$ at more than two points. Choose two points $p'$
and $p''$ among these intersection point. Let $v$ be the point at
which two lines $\overline{p_{12}, p'}$ and $\overline{p_{13},
p''}$ intersect. Then the cone over the curve $C$ with vertex $v$
has all the point of $\Gamma$ but not the point $q$.

\emph{Case 10.} (r=12)

All the points except one point, $p_{13}$, lie on the hyperplane
$H$. It immediately follows from Theorem~\xref{theorem-of-Bese}
that we can find a plane quintic curve  which passing
$\{p_1,\cdots, p_{12}\}$ but not the point $q$. Taking the cone
over the plane quintic curve with vertex $p_{13}$, we obtain a
quintic hypersurface that we want.

\emph{Case 11.} (r=13)

In this case, all the points lie on the hyperplane $H$. It
immediately follows from Theorem~\xref{theorem-of-Bese} that we
can find a plane quintic passing all the point except the point
$q$, which gives us a quintic hypersurface  in $\PP^3$ that we
want.

Consequently, we complete the proof.
\end{proof}

Therefore, the first part of Theorem~B has been proved.

The three-dimensional conjecture of Fujita (see \cite{EinLa93},
\cite{Ka97}, and \cite{Rei88}) implies
Lemma~\xref{lemma-on-blow-ups-of-points} in the case when the
points in $\Gamma$ are in \emph{very general} position. Moreover,
in the case when points in $\Sigma$ are in \emph{very general}
position, the $\QQ$-factoriality of $X$ follows from Lefschetz
theory (see Theorem~1.34 in \cite{Cl83}). However,
 neither three-dimensional nor two-dimensional conjecture of
Fujita cannot be, in general, applied to an appropriate adjoint
linear system in our case. The crucial point here is that the
proof of Theorem~\xref{theorem-of-Bese} is based on the vanishing
theorem of Ramanujam (see \cite{Bom73} and \cite{Ram72}) for
$2$-connected effective divisors on an algebraic surface (see
Proposition~2 in \cite{vdVe79}) which is slightly stronger in some
cases than the vanishing theorem of Kawamata and Viehweg (see
\cite{Ka82} and \cite{Vi82}).

The method of \cite{Cl83} also explains the non-$\QQ$-factoriality
of Examples~\ref{example-on-Q-factoriality},
~\ref{example-on-quartic-with-double-point},
and~\ref{example-on-bidegree-2-3}
 over $\CC$.
 Let
$X\longrightarrow \PP^3$ be a double cover ramified along a sextic
$S$. Suppose that the sextic $S\subset\mathbb{P}^3$ is given by
the equation
\[g^2_3(x,y,z,w)+h_r(x,y,z,w)f_{6-r}(x,y,z,w)=0,\]
where $g_3$, $h_r$, and $f_{6-r}$, $1\leq r\leq 3$, are generic
homogeneous polynomials over $\CC$ of degree $3$, $r$, and $ 6-r$,
respectively. Then the number of singular points is $18r-3r^2$.
All  of them are simple double points. The defect of $V$ is
\begin{equation*}\begin{split}h^1(\PP^3,
\mathcal{I}\otimes\mathcal{O}_{\PP^3}(5)) &= h^0(\PP^3,
\mathcal{I}\otimes\mathcal{O}_{\PP^3}(5))-
h^0(\PP^3,\mathcal{O}_{\PP^3}(5))+ h^0(\PP^3,\bigoplus_{p\in\Sigma}\CC)\\
&=h^0(\PP^3, \mathcal{I}\otimes\mathcal{O}_{\PP^3}(5))-
56+18r-3r^2.\\
\end{split}\end{equation*}
Let $H$ be the hypersurface of degree $r$ defined by $h_r=0$. Then
it is easy to check
\begin{equation*}\begin{split}
h^0(\PP^3,\mathcal{I}\otimes\mathcal{O}_{\PP^3}(5)) &\geq
h^0(\PP^3,\mathcal{O}_{\PP^3}(4))+h^0(H,\mathcal{O}_H(2))+h^0(H,\mathcal{O}_H)=42
 \mbox{ when $r=1$}\\
&\geq
h^0(\PP^3,\mathcal{O}_{\PP^3}(3))+h^0(H,\mathcal{O}_H(2))+h^0(H,\mathcal{O}_H(1))=33
\mbox{ when $r=2$}\\
&\geq
h^0(\PP^3,\mathcal{O}_{\PP^3}(2))+h^0(H,\mathcal{O}_H(2))+h^0(H,\mathcal{O}_H(2))=30
\mbox{ when $r=3$}.\\
\end{split}\end{equation*}
In all the cases, the defect of $V$ is positive. Therefore, the
double cover $X$ is not $\QQ$-factorial.


\section{Elliptic fibrations.}\label{section-of-elliptic}

This section is devoted to Theorem~C.

Let $\pi:X\to \PP^3$ be a $\QQ$-factorial double cover ramified
along a nodal sextic $S\subset \PP^3$. Consider a fibration
$\tau:Y\to Z$ whose general enough fiber is a smooth elliptic
curve. Suppose that we have a birational map $\rho$ of $X$ onto
$Y$. We then put $\MMM_{X}={\frac {1} {n}}\MMM$ with
$\MMM=\rho^{-1}(|\tau^{*}(H_Z)|)$, where $H_Z$ is a very ample
divisor on $Z$ and $n$ is the natural number such that
$\MMM\subset |-nK_{X}|$.
\[ \xymatrix{
X\ar[d]_{\pi}\ar@{-->}[rr]^{\rho}&& Y\ar[d]^{\tau}\\
\PP^3&&Z}\]%

It immediately follows from
Theorem~\xref{theorem-on-elliptic-fibrations} that the set
$\CCSS(X, \MMM_{X})$ is non-empty.

\begin{remark}
\label{remark-on-pencils} The linear system $\MMM$ is not composed
from a pencil and cannot be contained in the fibers of any
dominant rational map $\chi:X\dasharrow \PP^1$.
\end{remark}

Using the proof of Lemma~\xref{lemma-on-smooth-points}, we can
easily show that the set $\CCSS(X, \MMM_{X})$ does not contain any
smooth point of $X$.

\begin{lemma}
\label{lemma-on-singular-point} Let $o$ be a simple double point
on $X$ that belongs to $\CCSS(X, \MMM_{X})$. Then there is a
birational map $\beta:\PP^{2}\dasharrow Z$ such that the diagram
\[ \xymatrix{
X\ar[d]_{\pi}\ar@{-->}[rr]^{\rho}&& Y\ar[d]^{\tau}\\
\PP^3\ar@{-->}[r]^{\gamma}&\PP^2\ar@{-->}[r]^{\beta}&Z}\]%
commutes, where $\gamma$ is the projection from the point
$\pi(o)$.
\end{lemma}

\begin{proof}
Let $f:W\to X$ be the blow up at the point $o$ and $C$ be a
general enough fiber of the elliptic fibration
$\phi_{|-K_{W}|}:W\to \PP^{2}$. Then for a general surface $D$ in
$f^{-1}(\MMM)$,
$$
2(n-\mult_{o}(\MMM))=C\cdot D\ge 0,
$$
while $\mult_{o}(\MMM)\ge n$ by
Theorem~\xref{theorem-about-double-point}. We can therefore
conclude that  $\mult_{o}(\MMM)=n$ and $f^{-1}(\MMM)$ lies in the
fibers of the elliptic fibration $\phi_{|-K_{W}|}:W\to \PP^{2}$,
which implies the claim.
\end{proof}

\begin{corollary}
\label{corollary-about-two-singular-points} The set $\CCSS(X,
\MMM_{X})$ cannot contain two singular points of the $3$-fold $X$.
\end{corollary}

We  assume that $\CCSS(X, \MMM_{X})$ does not contain any point
and that it contains a curve $C\subset X$.

\begin{lemma}
\label{lemma-on-degree-of-curve-two} The intersection number
$-K_{X}\cdot C$ is $1$.
\end{lemma}

\begin{proof}
Let $H$ be a general enough divisor in the linear system
$|-K_{X}|$. Then we have
$$
2=H\cdot K_{X}^{2}= H\cdot \MMM_{X}^{2}\geq
\mult_{C}(\MMM_{X}^{2})H\cdot C\geq -K_{X}\cdot C,
$$
which implies $-K_{X}\cdot C\le 2$.

Suppose $-K_{X}\cdot C=2$. Then $\Supp(\MMM_{X}^{2})=C$ and
$$
\mult_{C}(\MMM_{X}^{2})=mult^2_{C}(\MMM_{X})=1,
$$
which means that for two different divisors $M_{1}$ and $M_{2}$ in
the linear system $\MMM$ we have
$$
\mult_{C}(M_{1}\cdot M_{2})=n^2,
\quad\mult_{C}(M_{1})=\mult_{C}(M_{2})=n,
$$
and set-theoretically $M_{1}\cap M_{2}=C$. However, the linear
system $\MMM$ is not composed from a pencil. Therefore, for a
general enough point $p\not\in C$ the linear subsystem $\DDD$ of
$\MMM$ passing through the point $p$ has no base components. Let
$D_{1}$ and $D_{2}$ be general enough divisors in $\DDD$. Then in
set-theoretic sense
$$p\in D_{1}\cap D_{1}=M_{1}\cap M_{2}=C,$$
which is a contradiction.
\end{proof}

\begin{corollary}
\label{corollary-on-image-of-curve-two} The curve $\pi(C)\subset
\PP^3$ is a line and $C\cong \PP^1$.
\end{corollary}

\begin{remark}
\label{remark-on-reducible-curves} In the second part of the proof
of Lemma~\xref{lemma-on-degree-of-curve-two}, we have never used
the irreducibility of the curve $C$. Hence, we may assume
$\CCSS(X, \MMM_{X})=\{C\}$. Moreover, the same arguments imply
$\mult_{C}(\MMM^2)<2n^{2}$.
\end{remark}

\begin{lemma}
\label{lemma-on-good-curve-two} The line $\pi(C)$ is contained in
the sextic surface $S$.
\end{lemma}

\begin{proof}
It follows from the proof of  Lemma~\ref{lemma-on-good-curve} and
Remark~\ref{remark-on-reducible-curves}
\end{proof}

Before we proceed, we observe
$$
\#|\Sing(X)\cap C|\leq\left\{\aligned
&3,\ \pi(C)\not\subset S\\
&5,\ \pi(C)\subset S,\\
\endaligned
\right.
$$
by intersecting $S$ with either the line $\pi(C)$ or a hyperplane
in $\PP^3$ passing through $\pi(C)$.  Furthermore, when
$\pi(C)\subset S$, the equality $\#|\Sing(X)\cap C|=5$ holds if
and only if all the hyperplanes tangent to the sextic surface $S$
at points of $\pi(C\backslash\Sing(X))$ coincide.

\begin{lemma}
\label{lemma-on-three-points-on-bad-curve} The curve $C$ passes
through at least four singular points of $X$.
\end{lemma}

\begin{proof}
Let $H$ be a general hyperplane in $\PP^3$ containing the line
$\pi(C)$. Then the curve
$$
D=H\cap S=\pi(C)\cup Q
$$
is reduced, where $Q$ is a quintic curve. The curve $D$ is
singular at each singular point $p_{i}$ of $S$ such that $p_{i}\in
\pi(C)$ for $i\in \{1,\ldots,k\}$. The set $\pi(C)\cap Q$ consists
of at most $5$ points and  $\Sing(D)\cap \pi(C)\subset \pi(C)\cap
Q$. Thus $k=\#|\Sing(X)\cap C|\leq 5$.

Suppose $k\le 3$. Then the intersection $\pi(C)\cap Q$ contains
two points $o_{1}$ and $o_{2}$ different from $p_{i}$ due to the
generality in the choice of $H$. The hyperplane $H$ is therefore
tangent to the sextic $S$ at $o_{1}$ and $o_{2}$. Hyperplanes
passing through the line $\pi(C)$ form a pencil whose proper
transforms on the $3$-fold $X$ are  K3 surfaces in $|-K_{X}|$
passing through $C$. Hence, the lines  tangent to the sextic
surface $S$ at a general point of the line $\pi(C)$ span whole
$\PP^{3}$. Note that this is no longer true in the case $k=5$ as
we mentioned right before the lemma.

Let $L_1$ and $L_2$ be general enough lines in $H$ passing through
the points $o_{1}$ and $o_2$, respectively. Then $L_{j}$ is
tangent to the sextic surface $S$ at the point $o_{j}$. Therefore,
the proper transform ${\tilde L}_{j}\subset X$ of the curve
$L_{j}$ is an irreducible curve such that $-K_{X}\cdot {\tilde
L}_{j}=2$. Also, it is singular at the point ${\tilde
o}_{j}=\pi^{-1}(o_{j})$. Consider the proper transform ${\tilde
H}$ of the surface $H$ on $X$ and a general surface $M$ in the
linear system $\MMM$. Then
$$
M\vert_{\tilde H}=\mult_{C}(\MMM)C+R,
$$
where $R$ is an effective divisor on ${\tilde H}$ such that
$C\not\subset\Supp(R)$. Moreover,
$$
2n=M\cdot {\tilde L}_{j}\ge \mult_{{\tilde o}_{j}}({\tilde
L}_{j})\mult_{C}(M)+\sum_{p\in (M\backslash C)\cap {\tilde L}_{j}}
\mult_{p}(M)\cdot \mult_{p}({\tilde L}_{j})\ge 2n,
$$
which implies $M\cap {\tilde L}_{j}\subset C$ set-theoretically.
However, on ${\tilde H}$ the curves ${\tilde L}_{1}$ and ${\tilde
L}_{2}$ span two pencils with the base loci consisting of the
points ${\tilde o}_{1}$ and ${\tilde o}_{2}$, respectively.
Therefore, we see $R=\emptyset$ due to the generality in the
choice of two curves $ L_{1}$ and $L_{2}$. Note that if $k=4$,
then this is not true.

Hence,  set-theoretically $M\cap {\tilde H}=C$ for a general
divisor ${\tilde H}\in |-K_{X}|$ passing through the curve $C$ and
a divisor $M\in\MMM$ with ${\tilde H}\not\subset\Supp(M)$. Let
${\tilde p}$ be a general point on the surface ${\tilde H}$ and
$\MMM_{{\tilde p}}$ be the linear system of surfaces in $\MMM$
containing ${\tilde p}$. Then $\MMM_{{\tilde p}}$ has no base
components due to Remark~\xref{remark-on-pencils}. Therefore, for
a general divisor ${\tilde M}$ in $\MMM_{\tilde p}$
$$
{\tilde p}\in {\tilde M}\cap {\tilde H}=C
$$
because ${\tilde H}\not\subset \Supp({\tilde M})$, which
contradicts the generality of the point ${\tilde p}\in {\tilde
H}$.
\end{proof}

\begin{lemma}
\label{lemma-on-four-points} Suppose that the curve $C$ contains
exactly $4$ singular points of the $3$-fold $X$. Then there is a
birational map $\beta:\PP^2\dasharrow Z$ such that the diagram
\[ \xymatrix{
X\ar@{-->}[d]_{\Xi_{\pi(C)}}\ar@{-->}[rr]^{\rho}&& Y\ar[d]^{\tau}\\
\PP^2\ar@{-->}[rr]^{\beta}&&Z}\] %
is commutative, where $\Xi_{\pi(C)}$ is a rational map defined as
in Construction~B.
\end{lemma}

\begin{proof}
By our assumption, the curve $C$ passes through  four singular
points $p_1$, $p_2$, $p_3$, $ p_4$ of $X$. We consider the blow up
$g_1:\tilde{X}\to X$ at the points $p_1, \cdots, p_4$ and the blow
up $g_2:W\to\tilde{X}$ along the proper transform of the curve $C$
on $\tilde{X}$. Put $g:=g_2\circ g_1:W\to X$. We then get
$$
-K_{W}=g^{*}(-K_{X})-\sum_{i=1}^{4}E_{i}-F,
$$
where $E_{i}$ and $F$ are the $g$-exceptional divisors such that
$g(E_{i})=p_i$ and $g(F)=C$. Let $L$ be a curve on $W$ such that
$\pi\circ g(L)$ is a line tangent to $S$ at some general point of
$\pi(C)$. Then
$$
\MMM_{W}\cdot L\le 2-2\mult_{C}(\MMM_{X})\le 0,
$$
where $\MMM_{W}=g^{-1}(\MMM_{X})$. Because such curves as $L$ span
a Zariski dense subset in $W$, we obtain $\mult_{C}(\MMM_{X})=1$.
Each elliptic curve $L$ is a fiber of the elliptic fibration
$\Xi_{\pi(C)}\circ g:W\to \PP^2$. Thus $\MMM_{W}$ lies in the
fibers of $\Xi_{\pi(C)}\circ g$, which implies the claim.
\end{proof}

\begin{lemma}
\label{lemma-on-five-points} The curve $C$ passes through at most
$4$ singular points of $X$.
\end{lemma}

\begin{proof}
Suppose that  the curve $C$ passes through  $5$ singular points
$p_1, \cdots, p_5$ of $X$. Again, we consider the blow up
$g_1:\tilde{X}\to X$ at the points $p_1, \cdots, p_5$ and the blow
up $g_2:W\to\tilde{X}$ along the proper transform of the curve $C$
on $\tilde{X}$. Put $g:=g_2\circ g_1:W\to X$. Then we obtain
$$
-K_{W}=g^{*}(-K_{X})-\sum_{i=1}^{5}E_{i}-F,
$$
where $E_{i}$ and $F$ are the $g$-exceptional divisors such that
$g(E_{i})=p_i$ and $g(F)=C$. Let $f:U\to W$ be a birational
morphism such that $h=\rho\circ g\circ f$ is a morphism. Then we
obtain
$$
K_{U}+\MMM_{U}=(g\circ f)^{*}(K_{X}+\MMM_{X})+\sum_{i=0}^{r}a_{i}G_{i},%
$$%
 where $\MMM_{U}=(g\circ f)^{-1}(\MMM_{X})$, $G_{i}$ are the $(g\circ f)$-exceptional
divisors, and $a_{i}\in \QQ$. Whenever $a_{i}\le 0$, we have
$g\circ f(G_{i})=C$. But $\mult_{C}(\MMM_{X})<2$ by
Remark~\xref{remark-on-reducible-curves} and hence there is
exactly one $i$, say $i=0$, such that $a_{0}\le 0$. It implies
$f(G_{0})=F$ and $a_{0}=0$.

Consider a general enough fiber ${\hat L}$ of the morphism
$\tau\circ h:U\to Z$. Then $K_{U}\cdot {\hat L}=0$ because the
curve ${\hat L}$ is elliptic. However, $\MMM_{U}\cdot {\hat L}=0$
by construction. So we see $G_{i}\cdot {\hat L}=0$ for $i\ne 0$,
which means that $f$ is an isomorphism near ${\hat L}$. Thus
$\MMM_{W}\cdot {\tilde L}=0$, where $\MMM_{W}= f^{-1}(\MMM_{X})$
and ${\tilde L}=f({\hat L})$.

There is a surface $D\subset W$ such that $\pi\circ
g(D)\subset\PP^3$ is the plane tangent to the sextic surface $S$
along the whole line $\pi(C)$. The surface $D$ is the closure of
the set spanned by curves whose images via $\pi\circ g$ are lines
tangent to the surface $S$ at some point of $\pi(C)$.

By the same argument as in the proof of
Lemma~\ref{lemma-on-four-points}, we obtain
$\mult_{C}(\MMM_{X})=1$, and hence
$$
D\sim \MMM_{W}-F+\sum_{i=1}^{5}b_{i}E_{i}
$$
for some $b_{i}\in\ZZ$.  On the other hand, because ${\hat L}\cdot
G_{i}=0$ for $i\ne 0$, we get
$$E_{j}\cdot {\tilde
L}= f^{*}(E_{j})\cdot{\hat L}=\sum_{i=1}^{r}c_{ij}G_{i}\cdot{\hat
L}=0
$$
where $c_{ij}\in \NN$. Therefore, ${\tilde L}\cdot D<0$, which
means ${\tilde L}\subset D$. This is impossible because the curves
${\tilde L}$ span a Zariski dense subset in $W$.
\end{proof}

Therefore, Theorem~C is proven.


\section{Canonical Fano $3$-folds.}

To prove Theorem~D, we let  $\pi:X\to \PP^3$ be a $\QQ$-factorial
double cover ramified in a nodal sextic $S\subset \PP^3$. We then
suppose that there is a non-biregular birational map
$\rho:X\dasharrow Y$ of $X$ onto a Fano $3$-fold $Y$ with
canonical singularities. We are to show that there is a curve
$C\subset X$ such that $\pi(C)$ is a line on the surface $S$
passing through five nodes of the sextic $S$.

We put $\MMM=\rho^{-1}(|-nK_{Y}|)$ and $\MMM_{X}={\frac {1}
{n}}\MMM$ for a natural number $n\gg 0$. We then see that
$K_{X}+\MMM_{X}\sim_\QQ 0$ and the singularities of the movable
log pair $(X, \MMM_{X})$ are not terminal by
Theorem~\xref{theorem-on-canonical-Fano-varieties}. By our
construction, the linear system $\MMM$ cannot be contained in the
fibers of any dominant rational map $\chi:X\dasharrow Z$ with
$0<\dim(Z)\leq 2$.

\begin{proposition}
\label{corollary-on-five-points}The set $ \CCSS(X, \MMM_{X})$
consists of a single curve $C\subset X$ which satisfies
\begin{enumerate}
\item$-K_{X}\cdot C=1$, \item$\pi(C)\subset S$, \item
 $\#|\Sing(X)\cap C|=5$.
 \end{enumerate}
\end{proposition}
\begin{proof}
For the proof, we  literally repeat the proofs in
Section~\ref{section-of-elliptic} except those of
Lemmas~\xref{lemma-on-singular-point} and
\xref{lemma-on-four-points}.
\end{proof}

Let $p_1, p_2, p_3, p_4, p_5\in C$ be  singular points of $X$. We
consider the blow up $f_1:\tilde{X}\to X$ at all the points
$p_{i}$ and the blow up $f_2:W\to\tilde{X}$ along the proper
transform of the curve $C$ on $\tilde{X}$. Put $f=f_2\circ
f_1:W\to X$. We then note that $W$ is smooth and
$$
-K_{W}\sim f^{*}(-K_{X})-\sum_{i=1}^{5}E_{i}-G,
$$
where $E_{i}$ and $G$ are the $f$-exceptional divisors with
$f(E_{i})=p_{i}$ and $f(G)=C$. Each surface $E_{i}$ is isomorphic
to the  blow up of $\PP^{1}\times \PP^{1}$ at one point.
We have the proper transforms $F^{i}_{1}$ and $F^{i}_{2}$  of two
rulings of the quadric $\PP^{1}\times \PP^{1}$ with
self-intersection $-1$ on each surface $E_i$.

The normal bundle $\NNN_{W/F^{i}_{j}}$ of the curve
$F^{i}_{j}\cong \PP^1$ in the $3$-fold $W$ splits into
$$
\NNN_{W/F^{i}_{j}}\cong\OOO_{F^{i}_{j}}(a)\oplus\OOO_{F^{i}_{j}}(b)%
$$
for some integers $a\geq b$. The exact sequence
$$
0\to \TTT_{F^{i}_{j}}\to \TTT_{W}\vert_{F^{i}_{j}}\to
\NNN_{W/F^{i}_{j}}\to 0
$$
implies $\deg(\NNN_{W/F^{i}_{j}})=a+b=-K_{W}\cdot
F^{i}_{j}+2g(F^i_j)-2=-2$. On the other hand, the exact sequence
$$
0\to \NNN_{E_{i}/F^{i}_{j}}\to \NNN_{W/F^{i}_{j}}\to \NNN_{W/E_{i}}\vert_{F^i_j}\to 0%
$$
together with  $\NNN_{E_{i}/F^{i}_{j}}\cong \OOO_{F_j^i}(-1)$
implies $b\geq -1$. Therefore, $a=b=-1$ and we can make a standard
flop for each curve $F^{i}_{j}$. Indeed, we let
$h:\widetilde{W}\to W$ be the blow up along all the curves
$F^{i}_{j}$ and $R^{i}_{j}$ be the $h$-exceptional divisor
dominating the curve $F^{i}_{j}$. Then $R^{i}_{j}\cong
\PP^{1}\times \PP^{1}$ and there is a birational morphism ${\hat
h}:\widetilde{W}\to {\hat W}$ which contracts each surface
$R^{i}_{j}$ to a curve ${\hat F}^{i}_{j}\subset {\hat W}$ and for
which ${\hat h}\circ h^{-1}$ is not an isomorphism in a
neighborhood of each curve $F^{i}_{j}$.

Let ${\hat E}_{i}={\hat h}\circ h^{-1}(E_{i})\subset {\hat W}$.
Then ${\hat E}_{i}\cong \PP^2$ and
$$
{\hat E}_{i}\vert_{{\hat E}_{i}}\cong \OOO_{\PP^2}(-2),
$$
which implies that each divisor ${\hat E}_{i}$ can be contracted
to a terminal cyclic quotient singularity of type ${\frac
{1}{2}}(1,1,1)$. Let ${\hat f}:{\hat W}\to V$ be the contraction
of all the ${\hat E}_{i}$. Then $V$ has exactly five singular
points $o_{i}$ of type ${\frac {1}{2}}(1,1,1)$, it is
$\QQ$-factorial, and $\Pic(V)\cong\ZZ\otimes \ZZ$.

Let $F={\hat f}\circ{\hat h}\circ h^{-1}(G)$. Then there is a
birational morphism $g:V\to X$ contracting the surface $F$ to the
curve $C$.
\[ \xymatrix{&\widetilde{W}\ar[dl]_{h}\ar[dr]^{\hat{h}}&\\
W\ar[d]_{f}\ar@{-->}[rr]^{\hat{h}\circ h^{-1}}&& \hat{W}\ar[d]^{\hat{f}}\\
X\ar@{-->}[d]_{\rho}&&V\ar[ll]_{g}\ar[d]^{\phi_{|-rK_V|}}\\
Y\ar@{-->}[rr]&&U}\]%
At a generic point of $C$ the morphism $g$ is a blow up. In fact,
the morphism $g$ is the blow up of the ideal sheaf of the curve
$C\subset X$ by Proposition 1.2 in \cite{Tz03}. Moreover, the
proof of Lemma~\xref{lemma-on-bad-curve} implies
$\mult_{C}(\MMM_{X})=1$. Hence,
$$
-K_{V}\sim_\QQ \MMM_{V}\sim_\QQ g^{*}(-K_{X})-F,
$$
where $\MMM_{V}=g^{-1}(\MMM_{X})$. The morphism $g\vert_{F}:F\to
C$ has five reducible fibers consisting of two copies of $\PP^{1}$
intersecting transversally at the corresponding singular point
$o_{i}$ that is a simple double point on the surface $F$.

Let ${\tilde C}\subset F$ be the unique base curve of the pencil
$|-K_{V}|$. Then the numerical equivalence ${\tilde C}\equiv
K_{V}^{2}$ holds. Therefore, we have
$$
-K_{V}\ {\text {is nef}}\ \iff -K_{V}\cdot {\tilde C}\geq 0\iff
-K_{V}^{3}\ge 0.
$$
Because elementary calculations imply $-K_{V}^{3}={\frac {1}{2}}$,
the anticanonical divisor $-K_{V}$ is nef and  big. Hence,
$|-rK_{V}|$ is base-point-free for a natural number $r\gg 0$ by
Base Point Freeness theorem (see \cite{KMM}). The morphism
$\phi_{|-rK_{V}|}:V\to U$ is birational and $U$ is a canonical
Fano $3$-fold with $-K_{U}^{3}={\frac {1}{2}}$.

The image of every element in the set $\CCSS(V, \MMM_{V})$ on the
$3$-fold $X$ is an element in $\CCSS(X, \MMM_{X})$ because
$K_{V}+\MMM_{V}=g^*(K_X+\MMM_X)$. Hence, every element in
$\CCSS(V, \MMM_{V})$ must be a curve dominating the curve $C$ due
to the assumption made in
Remark~\xref{remark-on-reducible-curves}, which implies
$\mult_{C}(\MMM)\ge 2n^{2}$. However, it is impossible because of
Remark~\xref{remark-on-reducible-curves}. Therefore, the set
$\CCSS(V, \MMM_{V})=\emptyset$.

For a rational number $c$ slightly bigger than $1$,  the
singularities of the log pair $(V, c \MMM_{V})$ are still terminal
and the equivalence
$$
K_{V}+c\MMM_{V}=
\phi_{|-rK_{V}|}^{*}(K_{U}+c\MMM_{U}), %
$$
holds, where $\MMM_{U}=\phi_{|-rK_{V}|}(\MMM_{V})$. Hence, the
movable log pair $(U, c \MMM_{U})$ is a canonical model. On the
other hand, the movable log pair $(Y, {\frac{c}{n}}|-nK_{Y}|)$ is
a canonical
model as well. Consequently, the map %
$$
\phi_{|-rK_{V}|}\circ (\rho\circ g)^{-1}:Y\dasharrow U
$$ %
is an isomorphism by
Proposition~\xref{uniqueness-of-canonical-model}.

 All the statements above do not depend on the existence of a birational map
$\rho$ of $X$ onto $Y$. They depend only on the condition that $X$
has a curve $C$ such that $\pi(C)\subset S$ is a line passing
through five nodes of the sextic surface $S$.
 Once such a curve $C\subset X$
exists, we can construct a birational transformation of $X$ into a
canonical Fano $3$-fold by means of blowing up the ideal sheaf of
the curve $C\subset X$ and the birational morphism given by a
plurianticanonical linear system.

We have proved Theorem~D. In addition, we have obtained explicit
classification of all birational transformations of a double cover
$X$ into Fano $3$-folds with canonical singularities.

 As we mentioned before, five singular points of the
surface $S$ lying on the line $\pi(C)\subset S$ force every
hyperplane in $\PP^3$ tangent to $S$ at some point of $\pi(C)$
smooth on $S$ to be tangent to the surface $S$ along  whole the
line $\pi(C)$. Such a tangent hyperplane is unique and its proper
transform on $V$ is the only divisor in the linear system
$|-K_{V}-F|$ which is contracted by the birational morphism
$\phi_{|-rK_{V}|}$ to a non-terminal point of the canonical Fano
$3$-fold $U$.


\section{Sextic double solids over finite fields.}

 We consider a double cover $\pi:X\to \PP^3$ defined over a perfect
field $\FF$ of characteristic $\cha(\FF)>5$. Suppose that the
$3$-fold $X$ is $\QQ$-factorial and that it is ramified along a
nodal  sextic surface $S\subset \PP^3$. Actually, we may assume
that the field $\FF$ is algebraically closed because $\FF$ is
perfect. We are to adjust the proofs of both Theorems~A and~C to
the case $\cha(\FF)>5$.

 We first list valid statements in Sections~\ref{section-birational-super-rigidity} and~\ref{section-of-elliptic}
 in the case $\cha(\FF)>5$.
 The following
remain valid:
\begin{enumerate}
\item Propositions~\xref{independence-of-Kodaira-dimension}, \xref{uniqueness-of-canonical-model}, and Theorem~\xref{theorem-of-Iskovskikh}; %
\item negativity of  exceptional loci (see \cite{Ar62} and Lemma 2.19 in \cite{Ko91}); %
\item resolution of singularities of $3$-folds (see \cite{Ab98} and \cite{Cos96}); %
\item numerical intersection theory on smooth $3$-folds (see \cite{Fu98}); %
\item elementary properties of blow ups (see \cite{Ha77}). %
\end{enumerate}

\begin{lemma}
\label{lemma-on-Nother-Fano-in-char-p}
Theorems~\xref{theorem-Nother-Fano-inequality} and
\xref{theorem-on-elliptic-fibrations} are valid in the case
$\cha(\FF)>5$.
\end{lemma}
\begin{proof}
The proofs for the case $\cha(\FF)=0$ depend only on the facts
listed above.
\end{proof}
 The following  may not
remain valid in the case $\cha(\FF)\ne 0$:
\begin{enumerate}
\item Theorem~\xref{theorem-about-double-point}; %
\item special cases of Bertini theorem (see \cite{GrHa78}). %

\end{enumerate}

For the birational super-rigidity, we need
Theorem~\ref{theorem-about-double-point} and Bertini theorem.

The characteristic-free  method for the proof of
Theorem~\ref{theorem-of-Iskovskikh} in \cite{Pu00} can be used to
prove Theorem~\xref{theorem-about-double-point}. However, we used
Theorem~\xref{theorem-about-double-point} just to prove
Lemmas~\xref{lemma-on-singular-points}. So instead of proving
Theorem~\xref{theorem-about-double-point} in the case
$\cha(\FF)>5$, we  prove Lemmas~\xref{lemma-on-singular-points}
only with Theorem~\ref{theorem-of-Iskovskikh}, which is enough for
the birational super-rigidity.

\begin{lemma}
\label{lemma-on-singular-points-in-char-p} Let $(X, \MMM_{X})$ be
a movable log pair such that $-(K_{X}+\MMM_{X})$ is ample and let
$o\in X$ be a simple double point.  Then the point $o$ does not
belong to $\CCSS(X,\MMM_{X})$.
\end{lemma}

\begin{proof}
Suppose  that the point $o$ belongs to the set
$\CCSS(X,\MMM_{X})$. Let $f:W\to X$ be the blow up at the point
$o$ and $C$ be a general enough fiber of the elliptic fibration
$\phi_{|-K_{W}|}:W\to \PP^{2}$. Then
$$
2(1-\mult_{o}(\MMM_{X}))>C\cdot \MMM_{W}\ge 0,
$$
where $\MMM_{W}=f^{-1}(\MMM_{X})$. This implies
$\mult_{o}(\MMM_{X})<1$.

We consider
$$K_{W}+\MMM_{W}=
f^{*}(K_{X}+\MMM_{X})+(1-\mult_{o}(\MMM_{X}))G,$$ where $G$ is the
$f$-exceptional divisor. We then see that there is a center $B\in
\CCSS(W, \MMM_{W})$ with $B\subset G$.

The intersection number of $\MMM_W$ with each ruling of $G\cong
\PP^1\times\PP^1$ is $\mult_o{\MMM_X}<1$. On the other hand, we
have $\mult_{B}(\MMM_{W})\ge 1$. Therefore, the center $B$ must be
a point and
$$
\mult_{B}(\MMM_{W}^{2})\ge 4
$$
by Theorem~\xref{theorem-of-Iskovskikh}.

Let $H_{1}$ and $H_{2}$ be two general surfaces in $|-K_{W}|$
passing through the point $B$. Then $H_{1}\cap H_{2}$ consists of
the fiber $E$ of the elliptic fibration $\phi_{|-K_{W}|}$ with
$B\in E$. Consider general enough divisors $D\in |-2K_{W}|$ and
$F_1, F_2\in |f^{*}(-K_{X})|$. Then the divisors $D$, $F_1$, and
$F_2$ do not pass through the point $B$ at all. The divisors
$H_{1}+F_{1}$, $H_{2}+F_{2}$, and $D+G$ are elements of the linear
subsystem $\HHH\subset |f^{*}(-2K_{X})-G|$ of surfaces passing
$B$. The intersection
$$
\Supp(H_{1}+F_{1})\cap \Supp(H_{2}+F_{2})\cap \Supp(D+G)
$$
contains $B$ and consists of finite number of points. In
particular, the linear system $\HHH$ has no base curves but $B$ is
a base point of $\HHH$. Let $H$ be a general surface in $\HHH$.
Then we obtain
$$4>H\cdot \MMM_{W}^{2}\ge \mult_{B}(H)\mult_{B}(\MMM_{W}^{2})\ge 4,$$
which is absurd.
\end{proof}

During excluding a one-dimensional member of $\CCSS(X,\MMM_{X})$,
we implicitly used Bertini theorem only one time just for the
following special case.

\begin{lemma}
\label{lemma-on-special-Bertini} Let $C\subset X$ be a curve with
$-K_{X}\cdot C=1$ and $\pi(C)\not\subset S$. Then a general enough
surface $H\in |-K_{X}|$ passing through $C$ is smooth along
$C\setminus\Sing(X)$.
\end{lemma}

\begin{proof}
The simple double points of the $3$-fold $X$ correspond to  the
simple double points of the sextic surface $S$ because
$\cha(\FF)\ne 2$. Meanwhile, the curve $L:=\pi(C)$  on $\PP^3$ is
a line. The line $L$ cannot pass through more than 3 singular
points of $S$; otherwise it would be contained in $S$.
The surface $D=\pi(H)\subset\PP^3$ is a plane containing $L$. The
singularities of surface $H$ correspond to the singularities of
the curve $D\cap S$ which is the ramification divisor of the
double cover $\pi:H\to D$. For a general enough surface $H\in
|-K_X|$, the plane $D$ is not tangent to the sextic $S$ at the
points of $L\setminus \Sing(S)$, which implies the claim.
\end{proof}

Therefore, the birational super-rigidity remains true over the
field $\FF$.

Now, we consider the statements in
Section~\ref{section-of-elliptic} over the field $\FF$. They also
require both Theorem~\ref{theorem-about-double-point} and Bertini
Theorem.

The reason why Theorem~\ref{theorem-about-double-point} is
required again is the lemma below. It can be however proved only
with Theorem~\ref{theorem-of-Iskovskikh}.

\begin{lemma}
\label{lemma-on-singular-point-in-char-p} Let $\rho:X\dasharrow Y$
be a birational map and $\tau:Y\to Z$ be a fibration whose general
fiber is a smooth elliptic curve. Let $(X, \MMM_X)$ be the movable
log pair such that $\MMM:=\rho^{-1}(|\tau^{*}(H)|)$ and
$\MMM_X=\frac{1}{n}\MMM$, where $H$ is a very ample divisor on
surface $Z$ and  $n$ is the natural number such that $\MMM\subset
|-nK_{X}|$. Suppose that the set $\CCSS(X,\MMM_{X})$ contains a
singular point $o\in X$. Then there is a birational map
$\beta:\PP^{2}\dasharrow Z$ such that the diagram
\[ \xymatrix{
X\ar[d]_{\pi}\ar@{-->}[rr]^{\rho}&& Y\ar[d]^{\tau}\\
\PP^3\ar@{-->}[r]^{\gamma}&\PP^2\ar@{-->}[r]^{\beta}&Z}\]%
commutes, where $\gamma$ is the projection from the point
$\pi(o)$.
\end{lemma}

\begin{proof}
Consider the blow up $f:W\to X$ at the point $o$. Let $C$ be a
general fiber of $\phi_{|-K_{W}|}$. Then
$$
2n-2\mult_{o}(\MMM)=C\cdot f^{-1}(\MMM)\ge 0,
$$
which implies $\mult_{o}(\MMM_{X})\le 1$. Furthermore, the
multiplicity $\mult_{o}(\MMM_{X})$ cannot be less than 1. Indeed,
if $\mult_{o}(\MMM_{X})<1$, then the proof of
Lemma~\ref{lemma-on-singular-points-in-char-p} shows a
contradictory inequalities
$$4\ge H\cdot \MMM_{W}^{2}\ge \mult_{B}(H)\mult_{B}(\MMM_{W}^{2})>4,$$
where $\MMM_W=f^{-1}(\MMM_X)$, $B$ is a center of $\CCSS(W,
\MMM_{W})$, and  $H$ is a general surface in $|f^{*}(-2K_{X})-E|$
passing through $B$.

In the case $\mult_{o}(\MMM_{X})=1$,  the linear system
$f^{-1}(\MMM)$ lies in the fibers of the elliptic fibration
$\phi_{|-K_{W}|}:W\to \PP^{2}$, which implies the claim.
\end{proof}

Bertini Theorem is required again only for the following statement
that can be proved without using Bertini theorem.
\begin{lemma}
\label{lemma-on-special-Bertini-two} Let $C$ be a curve on $X$
such that $-K_{X}\cdot C=1$ and $\pi(C)\subset S$. Suppose that $
\#|\Sing(X)\cap C|\le 3 $. Then a general surface $H\in |-K_{X}|$
passing through curve $C$ has at least $2$ different simple double
points on the curve $C\subset X$ at which the $3$-fold $X$ is
smooth.
\end{lemma}

\begin{proof}
The surface $\pi(H)\subset \PP^3$ is a plane passing through the
line $L:=\pi(C)\subset S$. Therefore,
$$
\pi(H)\cap S=L\cup Q,
$$
where $Q$ is a plane quintic. Whenever $H$ moves in the pencil of
surfaces in $|-K_{X}|$ passing through $C$, the quintic $Q$ moves
in a pencil on $S$ whose base locus is $\Sing(S)\cap L$. It gives
a finite morphism $\gamma:L\to \PP^{1}$ of degree
$5-\#|\Sing(S)\cap L|$ such that in the outside of the set
$\Sing(S)\cap L$ the morphism $\gamma$ is ramified at the points
where $L\cup Q$ is not a normal crossing divisor on the plane
$\pi(H)$. These points correspond to non-simple double points of
the surface $H$ contained in the curve $C$ and different from
$\Sing(X)\cap C$. However, this morphism can not be ramified
everywhere because we assumed $\cha(\FF)>5$.
\end{proof}
\begin{corollary}\label{corollary-number-of-points-cha}
Lemma~\ref{lemma-on-three-points-on-bad-curve} remains true in the
case $\cha(\FF)>5$.
\end{corollary}
\begin{proof}
Apply Lemma~\ref{lemma-on-special-Bertini-two} to the proof of
Lemma~\ref{lemma-on-three-points-on-bad-curve}.
\end{proof}

Because the proofs of Lemmas~\ref{lemma-on-four-points} and
\ref{lemma-on-five-points} are characteristic-free, Theorem~E is
true.


\section{Potential density.}\label{section-potential-density}

Now, we  prove Theorem~F.

Consider a double cover $\pi:X\to \PP^3$ defined over a number
field $\FF$ and ramified along a nodal sextic surface $S\subset
\PP^3$. We suppose that $\Sing(X)\ne\emptyset$. We will show that
the set of rational points of the $3$-fold $X$ is potentially
dense, which means that there exists a finite extension $\KK$ of
the field $\FF$ such that the set of all $\KK$-rational points of
the $3$-fold $X$ is Zariski dense.

 The rationality and the
unirationality of the $3$-fold $X$ over the field $\overline{\QQ}$
would automatically imply potential density of rational points on
$X$. However, the $3$-fold $X$ is non-rational in general due to
Theorem~A and the unirationality of the $3$-fold $X$ is unknown.
Moreover, $X$ is expected to be non-unirational in general.
Actually, the degree of a rational dominant map from $\PP^{3}$ to
a double cover of $\PP^3$ ramified in a very generic smooth sextic
surface must be divisible by $2$ and $3$ due to \cite{Ko96} and
\cite{Ko00}.

The following result was proved in \cite{BoTsch99}:

\begin{theorem}
\label{theorem-about-double-covers} Let $\tau:D\to \PP^2$ be a
double cover defined over a number field $\FF$ and ramified along
a reduced sextic curve $R\subset \PP^2$. Suppose
$\Sing(D)\ne\emptyset$. Then the set of rational points on the
surface $D$ is  potentially dense if and only if the curve
$R\subset\PP^2$ is not a union of six lines intersecting at a
single point.
\end{theorem}

Actually, Theorem~\xref{theorem-about-double-covers} is a special
case of the following result in \cite{BoTsch00}.

\begin{theorem}
\label{theorem-about-elliptic-KKK-surfaces} Let $D$ be a K3
surface defined over a number field $\FF$ such that $D$ has either
a structure of an elliptic fibration or an infinite group of
automorphisms. Then the set of rational points on $D$ is
potentially dense.
\end{theorem}

Hence, taking Theorem~C into consideration, we see that Theorem~F
is a three-dimensional analogue of
Theorem~\xref{theorem-about-double-covers}.

When singularities of the sextic surface $S$ are worse than simple
double points but  are not too bad, the double cover $X$ tends to
be  more rational (see \cite{Ch97}).
 Thus
Theorem~F must be true for sextic surfaces with any singularities
possibly except cones over sextic curves.
If the sextic surface $S\subset\PP^3$ is a reduced union of six
hyperplanes passing through one line, the set of rational points
on $X$ is  not potentially dense due to Faltings Theorem
(\cite{Fa83} and \cite{FaWu84}) because the $3$-fold $X$ is
birationally isomorphic to a product $\PP^{2}\times C$, where $C$
is a smooth curve of genus $2$.

As a matter of fact, the sets of rational points are  potentially
dense on  double covers of $\PP^{n}$ ramified along general enough
sextic hypersurfaces for $n\gg 0$ due to the following result
(\cite{CoMaMu02}):

\begin{theorem}
\label{theorem-on-unirationality-of-double-spaces} Let $V$ be a
double cover of $\PP^{n}$ ramified in a sufficiently general
hypersurface of degree $2d>4$. Then $V$ is unirational if $n\ge
c(d)$, where $c(d)\in\NN$ depends only on $d$.
\end{theorem}

We will prove the potential density of the set of rational points
on $X$ using the technique of \cite{BoTsch98}, \cite{BoTsch99},
and \cite{HaTsch00} which relies on the following result proved in
\cite{Mer96}.

\begin{theorem}
\label{theorem-of-Merel} Let $\FF$ be a number field. Then there
is an integer $n_{\FF}$ such that no elliptic curve defined over
$\FF$ has a $\FF$-rational torsion point of order $n>n_\FF$.
\end{theorem}

Let $o$ be a simple double point on $X$. The point $\pi(o)$ is a
node of the sextic surface $S$. Replacing the field $\FF$ by a
finite extension of $\FF$, we may assume that the point $o$ and
some other finitely many points that we will need in the sequel
are defined over $\FF$. Let $f:V\to X$ be the blow up at the point
$o$ with $f$-exceptional divisor $E$. Then the linear system
$|-K_{V}|$ is free and the morphism
$$
\phi_{|-K_{V}|}:V\to \PP^{2}
$$
is an elliptic fibration. The surface $E$ is a multisection of
$\phi_{|-K_{V}|}$ of degree $2$. Let $H$ be a ge\-neral surface in
$|-f^{*}(K_{X})|$. Then $H$ is a multisection of $\phi_{|-K_{V}|}$
of degree $2$ as well.

The following lemma is a corollary of Proposition 2.4 in
\cite{BoTsch98}.

\begin{lemma}
\label{lemma-saliently-ramified-multisection} Suppose that there
is a multisection $M$ of $\phi_{|-K_{V}|}$ of degree $d\geq 2$
such that the morphism $\phi_{|-K_{V}|}\vert_{M}$ is branched at a
point $p\in M$ which is contained in a smooth fiber of the
elliptic fibration $\phi_{|-K_{V}|}$. Then the divisor
$p_{1}-p_{2}\in\Pic(C_{b})$ is not a torsion for some distinct two
points $p_{1}$ and $p_{2}$ of the intersection $M\cap C_{b}$,
where $C_{b}=\phi_{|-K_{V}|}^{-1}(b)$ and $b$ is a $\CC$-rational
point in the complement to a countable union of proper Zariski
closed subsets in $\PP^{2}$.
\end{lemma}

\begin{proof}
See \cite{BoTsch98}.
\end{proof}

\begin{lemma}
\label{lemma-on-potential-density-modulo-multisection} Let $M\in
|H|$ be an irreducible multisection of $\phi_{|-K_{V}|}$ of degree
$2$ defined over $\FF$ such that the set of rational points on $M$
is  potentially dense in $M$ and $\phi_{|-K_{V}|}\vert_{M}$ is
branched at a point contained in a smooth fiber of
$\phi_{|-K_{V}|}$. Then the set of rational points on $X$ is
potentially dense.
\end{lemma}

\begin{proof}
For each $n\in\NN$, we let $\Phi_n$ be the set of points $p$ of
$M$ satisfying the following two conditions:
\begin{enumerate}
\item the point $p$ is contained in a smooth fiber $C_p$ of the
elliptic fibration $\phi_{|-K_{V}|}$; \item $ 2np= nH\vert_{C_{p}}
$ in $\Pic(C_{p})$.
\end{enumerate}
Let ${\overline \Phi}_{n}$ be the Zariski closure of the set
$\Phi_{n}$ in $M$.

Suppose ${\overline \Phi}_{n}=M$ for some $n$. Take a very general
fiber $C$ of $\phi_{|-K_{V}|}$ and let
$$
C\cap M=\{p_{1}, p_{2}\},
$$
where $p_{1}\ne p_{2}$. Then either $2np_{1}\sim nH\vert_{C}$ or
$2np_{2}\sim nH\vert_{C}$ because ${\overline \Phi}_{n}=M$.
However, $p_{1}+p_{2}\sim H\vert_{C}$. Thus
$$
2np_{1}\sim 2np_{2}\sim nH\vert_{C}
$$
and the element $p_{1}-p_{2}$ is a torsion in $\Pic(C)$.
Therefore, the $\CC$-rational point $\phi_{|-K_{V}|}(C)$ is
contained in the countable union of proper Zariski closed subsets
in $\PP^{2}$ of
Lemma~\xref{lemma-saliently-ramified-multisection}, which
contradicts the very general choice of the fiber $C$. Accordingly,
the set $\Phi_{n}$ is not Zariski dense in $M$ for any $n\in \NN$.
Moreover, it follows from Theorem~\xref{theorem-of-Merel} that
each set $\Phi_{n}$ for $n>n_\FF$, where $n_\FF$ is the number
defined in Theorem~\xref{theorem-of-Merel}, is disjoint from the
set of $\FF$-rational points on $M$.

Because of the assumption on the multisection $M$, we may assume
that the set of $\FF$-rational points on the surface $M$ is
Zariski dense. Take an $\FF$-rational point
$$
q\in M':=M\backslash(Z\cup_{i=1}^{n_\FF} {\overline \Phi}_{i}),
$$
where the set $Z\subset M$ consists of points contained in
singular fibers of $\phi_{|-K_{V}|}$. Let $C_{q}$ be the fiber of
$\phi_{|-K_{V}|}$ passing through $q$. Then both the curve $C_{q}$
and the point $\phi_{|-K_{V}|}(q)$ are  defined over the field
$\FF$. The divisor $2q-H\vert_{C_{q}}\in \Pic(C_{q})$ is defined
over $\FF$ as well. Moreover, $2q-H\vert_{C_{q}}$ is not a torsion
divisor. By Riemann-Roch theorem, for each $n\in\NN$ there is a
unique $\FF$-rational point $q_{n}\in C_{q}$ such that
$$
q_{n}+(2n-1)q= nH\vert_{C_{q}}
$$
in $\Pic(C_{q})$. Because $2q-H\vert_{C_{q}}$ is not a torsion
divisor, we see that $q_{i}\ne q_{j}$ if and only if $i\ne j$. We
obtain an infinite collection of $\FF$-rational points on $C_q$.
Consequently, for each $\FF$-rational point $q$ in $M'$, the curve
$C_{q}$ is contained in the Zariski closure of the set of
$\FF$-rational points of $V$. Because the set $M'$ is a Zariski
dense subset of $M$, the set of rational points on the $3$-fold
$X$ is  potentially dense.
\end{proof}

In order to prove Theorem~F, it is enough to find an element in
$|H|$ satisfying the conditions of
Lemma~\ref{lemma-on-potential-density-modulo-multisection}. To
find such an element, we let $T$ be the set of points $(p,q)\in
S\times S$ satisfying the following conditions:
\begin{enumerate}
\item  $p\ne q$; %
\item the points $p$ and $q$ are smooth points on the sextic surface $S$;%
\item the point $q$ is contained in the hyperplane $D\subset\PP^3$ tangent to $S$ at  $p$; %
\item the point $q$ is a smooth point of the intersection $S\cap D$; %
\item the intersection $S\cap D$ is reduced. %
\end{enumerate}
We also let $\psi:T\to S$ be the projection on the second factor.
\begin{lemma}
\label{lemma-on-Gauss-map} The image $\psi(T)$ contains a Zariski
open subset of the sextic $S\subset\PP^3$.
\end{lemma}

\begin{proof}
Let $p$ be a general point on the sextic $S\subset\PP^3$ and $D$
be the hyperplane tangent to the surface $S$ at the point $p$ in
$\PP^3$. To prove the claim we just need to show that $D\cap S$ is
reduced, which is nothing but the finiteness of the Gauss map at a
generic point of $S$.

When the surface $S$ is smooth, the intersection $D\cap S$ is
known to be reduced (see \cite{FuLa81}, \cite{Ish82}, or
\cite{Pu95}). Even though $S$ can have  double points in our case,
the intersection $D\cap S$ is reduced because $S$ is not ruled
(see \cite{Mo77}).  Here, we prove it only with simple
calculation.

Suppose that $D\cap S$ is not reduced and
$$
D\cap S=mC+F\subset D\cong\PP^{2},
$$
where $m\ge 2$. Then $C$ is  a line, a conic, or a plane cubic
curve. Let $\gamma:{\tilde S}\to S$ be the blow up at the double
points of $S$ and ${\tilde C}=\gamma^{-1}(C)$. Then $S$ is a
surface of general type,
$$
K_{\tilde S}=\gamma^{*}(\OOO_{\PP^3}(2)\vert_{S}),
$$
and ${\tilde C}$ is either a rational curve or an elliptic curve.
Moreover, the self-intersection number ${\tilde C}^{2}$ of
$\tilde{C}$ is negative by adjunction formula, but ${\tilde C}$
moves in a family on the surface ${\tilde S}$ when we move the
point $p$ in $S$, which is a contradiction.
\end{proof}

Therefore, by Lemma~\xref{lemma-on-Gauss-map} we can find a
hyperplane $D\subset\PP^3$ such that $D\cap S$ is reduced and
singular at some smooth point of $S$. Moreover, we may assume that
$D$ does not contain the point $\pi(o)$ and there is a line
$L\subset \PP^{3}$ passing through the point $\pi(o)$ such that
$$
L\cap D\cap S\ne\emptyset
$$
and $L$ intersects the sextic $S$ transversally at four different
smooth points of $S$. Let ${\tilde D}$ be the surface in the
linear system $|H|$ such that $\pi\circ f({\tilde D})=D$. Then
${\tilde D}$ is an irreducible multisection of the elliptic
fibration $\phi_{|-K_{V}|}$ of degree $2$ such that
$\phi_{|-K_{V}|}\vert_{{\tilde D}}$ is branched at a point
$q\in{\tilde D}$ contained in the fiber $C$ of $\phi_{|-K_{V}|}$
such that $\pi\circ f(C)=L$. By construction, the fiber $C$ is a
smooth elliptic curve, $\pi\circ f(q)=L\cap D\cap S$, and $q$ is a
smooth point on ${\tilde D}$. Moreover, extending the field $\FF$
 we can assume that ${\tilde D}$ is defined over $\FF$. Hence, the
 set of
rational points is potentially dense on ${\tilde D}$ by
Theorem~\xref{theorem-about-double-covers}.  Theorem~F is proven.

It would be natural to prove Theorem~F in the case when the sextic
$S$ is singular and reduced (see
Theorem~\xref{theorem-about-double-covers}). Most of the arguments
in this section work for any reduced singular sextic surface.
Actually, in the case when the sextic $S$ has non-isolated
singularities (for example, when it is reducible) we do not need
to use Lemma~\xref{lemma-on-Gauss-map} at all, but in the case
when the sextic $S$ is irreducible and has isolated singularities
we can prove Lemma~\xref{lemma-on-Gauss-map} using the finiteness
of the Gauss map for curves (see \cite{Ha92}) in the assumption
$S$ is not a scroll (see \cite{Mo77}, \cite{Zak87}, and
\cite{Zak93}), which is satisfied automatically if $S$ is not a
cone. Moreover, in general the proof of Theorem~F must be simpler
for bad singularities. For instance, in the case when the sextic
$S$ has a singular point of multiplicity $4$, the double cover $X$
is unirational and non-rational in general due to \cite{Tu78}, but
it is rational when $S$ has a singular point of multiplicity $5$.
However, when $S$ is a cone over a smooth sextic curve
$R\subset\PP^{2}$, the double cover $X$ is birationally equivalent
to $\PP^{1}\times D$, where $D$ is a double cover of $\PP^{2}$
ramified along $R$. The potential density of rational points on
$X$ is therefore equivalent to  the potential density of rational
points on $D$, which is still unknown in general (see
\cite{BoTsch00}).


\end{document}